\documentclass{article}
   \usepackage{amssymb} 
 \usepackage[totalwidth=410pt, totalheight=600pt]{geometry}

\newcommand{\h}{\mbox{\textbf{h}}}

\newcommand{\al}{\alpha}
\newcommand{\eps}{\epsilon}
\newcommand{\la}{\lambda}
\newcommand{\La}{\Lambda}

\newcommand{\Th}{\Theta}

\newcommand{\mb}{\mbox}

\newcommand{\Mklz}[2]{\left\{\left. \, #1 \,\right|\, #2\, \right\}}
\newcommand{\W}{\mbox{$\Delta$}}

\newcommand{\pW}{\mbox{$\Delta^+$}}
\newcommand{\rW}{\mbox{$\Delta_{re}$}}
\newcommand{\nrW}{\mbox{$\Delta_{re}^-$}}
\newcommand{\prW}{\mbox{$\Delta_{re}^+$}}

\newcommand{\F}{\mathbb{F}}

\newcommand{\N}{\mathbb{N}}
\newcommand{\Nn}{\mathbb{N}_0}

\newcommand{\R}{\mathbb{R}}

\newcommand{\We}{\mbox{$\mathcal W$}}

\newcommand{\Cq}{\mbox{$\overline{C}$}}
\newcommand{\FTq}{\mbox{$\overline{F}_\Theta$}}
\newcommand{\RkX}{\mbox{$\mb{Fa}(X)$}}
\newcommand{\TD}{\mbox{$\widehat{T}$}}
\newcommand{\ND}{\mbox{$\widehat{N}$}}
\newcommand{\GD}{\mbox{$\widehat{G}$}}

\newcommand{\WeD}{\mbox{$\widehat{\We}$}}
\newcommand{\PD}{\mbox{$\widehat{P}$}}
\newcommand{\LD}{\mbox{$\widehat{L}$}}
\newcommand{\ve}[1]{\mbox{$\varepsilon (#1)$}} 
\newcommand{\ti}{\widetilde}

\newcommand{\red}[1]{red\left(#1\right)}
\newcommand{\goodone}{{\,\mbox{}_\triangleright\,}}
\newcommand{\goodtwo}{{\,\mbox{}_\triangleleft\,}}
\newcommand{\bad}{{\,\mbox{}_\Box\,}}

\newcommand{\Proof}{\mbox{\bf Proof: }}
\newcommand{\qed}{\mb{}\hfill\mb{$\square$}\\}

\newtheorem{thm}{Theorem}[section]
\newtheorem{defn}[thm]{Definition}
\newtheorem{prop}[thm]{Proposition}
\newtheorem{cor}[thm]{Corollary}
\newtheorem{rem}[thm]{Remark}

\newtheorem{lem}[thm]{Lemma}

\begin{document}

\title{Actions of the face monoid associated to a Kac-Moody group on its building}

\author{Claus Mokler\thanks{Supported by the Deutsche Forschungsgemeinschaft, and by the Marie Curie Research Training Network ``Liegrits'' of the European Community, MRTN-CT-2003-505078.}\\\\ Universit\"at Wuppertal, Fachbereich C - Mathematik\\  Gau\ss stra\ss e 20\\ D-42097 Wuppertal, Germany\vspace*{1ex}\\ 
          mokler@math.uni-wuppertal.de}
\date{}
\maketitle
\begin{abstract}\noindent 
We described in \cite{M1} a monoid $\GD$ acting on the integrable highest weight modules of a symmetrizable Kac-Moody algebra. It has similar structural properties 
as a reductive algebraic monoid with unit group a Kac-Moody group $G$. Now we find natural extensions of the action of the Kac-Moody group $G$ on its building 
$\Omega$ to actions of the monoid $\GD$ on $\Omega$. These extensions are partly motivated by representation theory and the combinatorics of the faces of the Tits cone.
\end{abstract}
\vspace*{1ex}\\
\section*{Introduction}
The face monoid $\GD$ is an infinite-dimensional algebraic monoid. It has been obtained in \cite{M1} by a Tannaka reconstruction from categories determined by the integrable highest weight representations of a symmetrizable Kac-Moody algebra. Compare also \cite{M7}, \cite{M8}. By its construction and by the involved categories it is a very natural object. Its Zariski open dense unit group $G$ coincides, up to a slightly extended maximal torus, with the special Kac-Moody group defined representation theoretically in \cite{KP1}.  

The face monoid is a purely infinite-dimensional phenomenon, quite unexpected. In the classical case, i.e., if you take a semisimple Lie algebra for the symmetrizable Kac-Moody algebra, it coincides with a semisimple simply connected algebraic group. Put in another way, there seem to exist fundamentally different infinite-dimensional generalizations of a semisimple simply connected algebraic group.

The results obtained in \cite{M1}, \cite{M2}, \cite{M3} and \cite{M4} show that the face monoid $\GD$ has similar structural and algebraic geometric properties as a reductive algebraic monoid, e.g., the monoid of ($n\times n$)-matrices. The face monoid $\GD$ is the first example of an infinite-dimensional reductive algebraic monoid. Actually, it is is particular. The investigation of the conjugagy classes in \cite{M5} will show, that the relation between the face monoid $\GD$ and its unit group $G$, the Kac-Moody group, is much closer than for a general reductive algebraic monoid.

Obviously, there is the following question: Does the face monoid $\GD$ fit in some way into the building geometric theory of the Kac-Moody group $G$?

In this article we investigate how to extend the action of the Kac-Moody group $G$ on its building $\Omega$ in a natural way to actions of the monoid $\GD$ on $\Omega$. 
For the article all the background from representation theory is not needed. We only use the algebraic properties of the monoid $\GD$, which in the end reduce to basic group-theoretical properties of the Kac-Moody group $G$ and properties of the faces of the Tits cone $X$. Note that we use the term "face" for a face of the convex cone $X$ in the sense of convex geometry. The structure of $\GD$ and the properties of the faces of the Tits cone $X$ are explained briefly to make the article easy to read.

To describe roughly how we extend the action of the Kac-Moody group $G$ on its building $\Omega$ to actions of the face monoid $\GD$ on $\Omega$, note the following property:
For the face monoid $\GD$ an infinite Renner monoid $\WeD$ plays the same role as the Weyl group $\We$ does for the Kac-Moody group $G$. For example, there are Bruhat and Birkhoff decompositions of $\GD$, similar as for $G$, but the Weyl group $\We$ replaced by the monoid $\WeD$. The monoid $\WeD$ can be easily constructed from the Weyl group $\We$ and the face lattice of the Tits cone $X$. In particular, the Weyl group $\We$ is the unit group of $\WeD$.

Now the building $\Omega$ associated to the Kac-Moody group $G$ is covered by certain subcomplexes, the apartments, which are isomorphic to the Coxeter complex $\mathcal C$ associated to the Weyl group $\We$. This connects the action of $G$ on $\Omega$ and the action of $\We$ on $\mathcal C$. We call an action of $\GD$ on $\Omega$, which extends the action of $G$, and an action of $\WeD$ on $\mathcal C$, which extends the action of $\We$, compatible if there holds a similar connection.

We show the following beautiful theorem: For an action of $\WeD$ on $\mathcal C$, extending the action of $\We$, satisfying certain necessary and sufficient conditions, there exists a uniquely determined compatible action of $\GD$ on $\Omega$. Every action of $\GD$ on $\Omega$, which extends the action of $G$, is obtained in this way. Furthermore, the action of $\GD$ preserves the natural order on $\Omega$ if and only if the action of $\WeD$ preserves the natural order on $\mathcal C$.

An action of $\WeD$ on $\mathcal C$, extending the action of $\We$, can be found in a natural way by representation theory.
Furthermore, there exists (a unique) compatible action of $\GD$ on $\Omega$. Unfortunately, these actions don't preserve the natural order on $\mathcal C$ resp. $\Omega$. For this reason we call these actions the bad actions for short.

We then look for modifications of the bad actions which avoid this defect. We call such actions good actions. Unexpectedly, we find not only one but two good actions of $\WeD$ on the Coxeter complex $\mathcal C$, for which there exist (uniquely determined) compatible good actions of $\GD$ on $\Omega$. To do this, we use the geometry of the Tits cone $X$, and the geometry of a modified Tits cone $X^\vee_{Loo}$ constructed by E. Looijenga. In a subsequent paper \cite{M6} we give an interesting algebraic geometric model for one of the good actions of $\GD$ on $\Omega$.

Some results which we obtain to carry out our investigations are important for itself. We give explicit descriptions of the lattice join and lattice intersection of the faces of the Tits cone. This leads to an elegant explicit description for the multiplication of two elements of the face monoid $\GD$. It also leads to a variant of a monoid theoretical construction of $\GD$, which is partly based on the root groups data system of $G$. All three things have been missing in \cite{M1}. Furthermore, we find Levi decompositions of the parabolic submonoids of $\GD$.

It remains open to determine all good actions of the face monoid $\GD$. Maybe these actions can be useful to find a building geometrical interpretation of $\GD$.

By its representation theoretic construction, the unit group of the face monoid $\GD$ is a Kac-Moody group associated to a 
special Kac-Moody root datum. The monoid theoretical construction of $\GD$ just mentioned may be very useful to find the correct definition of a face monoid whose unit group is a Kac-Moody group associated to a general Kac-Moody root datum.
It remains open to investigate and work out this case. 
%
%
%
%
%%%%%%%%%%%%%%%%%%%%%%%%%%%%%%%%%%%%%%%%%%%%%%%%%%%%%%%%%%%%%%%%%%%%%%%%%%%%%%%%%%%%%%%%%%%%%%%%%%%%%%%%%%%%%%%%%%%%%%%%%%%%%%%%%%%%%%%%%%%%%%%%%%%%%
%
%%%%%%%%%%%%%%%%%%%%%%%%%%%%%%%%%%%%%%%%%%%%%%%%%%%%%%%%%%%%%%%%%%%%%%%%%%%%%%%%%%%%%%%%%%%%%%%%%%%%%%%%%%%%%%%%%%%%%%%%%%%%%%%%%%%%%%%%%%%%%%%%%%%%%%
%
%
\section{Preliminaries} For our notation on Kac-Moody algebras and Kac-Moody groups we refer to \cite{M1}, Chapter 1, Preliminaries, or to \cite{M2}, Section 1, Preliminaries. We work here with a Kac-Moody group $G$ over a field $\F$ of characteristic zero whose derived group $G'$ coincides with the special Kac-Moody group defined and investigated by V. Kac and D. Peterson in \cite{KP1}. For our notation on convex cones and their faces we refer to \cite{M1}, Chapter 1, Preliminaries.
%
%
%
%%%%%%%%%%%%%%%%%%%%%%%%%%%%%%%%%%%%%%%%%%%%%%%%%%%%%%%%%%%%%%%%%%%%%%%%%%%%%%%%%%%%%%%%%%%%%%%%%%%%%%%%%%%%%%%%%%%%%%%%%%%%%%%%%%%%%%%%%%%%%%%%%% 
%
%
\section{The face lattice of the Tits cone\label{section2}}
Recall that $X=\We \overline{C}$ denotes the dual Tits cone in the real linear space $\h_\R^*$. It can be obtained by applying the Weyl group $\We$ to the closed fundamental chamber $\overline{C}:=\{\la\in\h_\R^*\mid\la(h_i)\geq 0 \mb{ for all }i\in I\}$. 
Set $C:= \{\la\in\h_\R^*\mid \la(h_i)> 0 \mb{ for all }i\in I\}$. For $J\subseteq I$ we call 
\begin{eqnarray*}
   F_J  &:=&   \{\la\in\h_\R^*\mid \la(h_i)=0 \mb{ for }i\in J,\;\la(h_i)>0 \mb{ for }i\in I\setminus J\},\\
   \overline{F}_J &:=& \{\la\in\h_\R^*\mid \la(h_i)=0 \mb{ for }i\in J,\;\la(h_i)\geq 0 \mb{ for }i\in I\setminus J \}
\end{eqnarray*}
the {\it open} resp. {\it closed facet} of $\overline{C}$ of type $J$.  By applying the elements of the Weyl group $\We$ to these facets, we get the {\it open} resp. {\it closed facets} of $X$ of type $J$. We do not call these sets "faces", which is more common in the literature, because we use "face" in a different meaning.

By a {\it face} of the Tits cone $X$ we mean a face of the convex cone $X$ in the sense of convex geometry. Denote the set of faces of the Tits cone $X$ by $\RkX$. Order $\RkX$ partially by the inclusion of faces, and equip it with the corresponding complete 
lattice structure. The action of the Weyl group $\We$ on the Tits cone $X$ induces an action of $\We$ on the partially ordered set of faces $\RkX$ resp. on the face lattice $\RkX$ of $X$. 

In this section we collect from the literature a description of the $\We$-action on $\RkX$, and a description of the partially ordered set $\RkX$, which are both important for this article. 
As new results, also important for the following sections, we obtain explicit descriptions of the lattice operations of $\RkX$. In addition, we describe  the intersections of the faces of the Tits cone with the closed facets of the Tits cone explicitely.
%
%%%%%%%%%%%%%%%%%%%%%%%%%%%%%%%%%%%%%%%%%%%%%%%%%%%%%%%%%%%%%%%%%%%%%%%%%%%%%%%%%%%%%%%%%%%%%%%%%%%%%%%%%%%%%%%%%%%%%%%%%%%%%%%%%%%
%Please add some additional space:
\vspace*{1ex}
%
%%%%%%%%%%%%%%%%%%%%%%%%%%%%%%%%%%%%%%%%%%%%%%%%%%%%%%%%%%%%%%%%%%%%%%%%%%%%%%%%%%%%%%%%%%%%%%%%%%%%%%%%%%%%%%%%%%%%%%%%%%%%%%%%%%%%
%

The following theorem and its corollary describe the set of faces $\RkX$ of $X$ as well as the $\We$-action on $\RkX$. These results are mainly a summary of some theorems of E. Looijenga in \cite{Loo} and P. Slodowy in \cite{Sl}. For detailed references please compare Section 2.1 of \cite{M1}.

To state the theorem and its corollary some definitions are useful. For $\emptyset\neq\Th \subseteq I$ we denote by $\Th^0$, respectively $\Th^\infty$ the set of indices which correspond to the sum of the components of the generalized Cartan submatrix $A_\Th$ of finite, respectively non-finite type. We set $\emptyset^0:=\emptyset^\infty:=\emptyset$. For a face $R$ of the Tits cone $X$ we denote by $ri(R)$ its relative interior. 

\begin{defn}\label{FL1}
A subset $\Th\subseteq I$ is called special if $\Th=\Th^\infty$.
A face $R$ of $X$ is called special if $ri(R)\cap \Cq\neq\emptyset$.
\end{defn}

For $J\subseteq I$ set $J^\bot:=\Mklz{i\in I}{a_{ij}=0\;\mb{ for all }j\in J}$. 
\begin{thm} \label{FL2} Let $\Th$ be a special subset of $I$. Then 
\begin{eqnarray}
   R(\Th):=\We_{\Th^\bot}\FTq
\end{eqnarray}
is a special face of $X$. Its relative interior is given by 
\begin{eqnarray}
 ri \left( R(\Th) \right) = \We_{\Th^\bot} \bigcup_{J\subseteq \Th^{\bot}\,,\, J=J^0} F_{\Th\cup J}.
\end{eqnarray}
The linear hull of $R(\Th)$ is given by
\begin{eqnarray}
  \qquad \qquad \mb{span}(R(\Th)) = \Mklz{\la\in \h_\R^*}{\la(h_i)=0\;\mb{ for all }\;i\in\Th}.
\end{eqnarray}
In particular, $R(\Th)=X\cap \mb{span}(R(\Th))=\Mklz{\la\in X}{\la(h_i)=0\;\mb{ for all }\;i\in\Th}$.\vspace*{1ex}

The $\We$-stabilizers of $R(\Th)$, pointwise and as a whole, are:
\begin{eqnarray}
      Z_{\mathcal W}(R(\Th)) &:=&\Mklz{\sigma\in\We}{\,\sigma\la=\la\mb{ for all }\la\in R(\Th)} = \We_{\Th}. \label{Stpoint}\\
      N_{\mathcal W}(R(\Th))&:=&\Mklz{\sigma\in\We}{\sigma R(\Th) = R(\Th)} = \We_{\Th\cup \Th^\bot}.\label{Stwhole}
\end{eqnarray}
\end{thm}

From this Theorem follows easily:
\begin{cor}\label{FL3}
(a) We get a bijective map from the special sets of $I$ to the special faces of $X$ by assigning to the special set $\Th$ the special face $R(\Th)$.\vspace*{1ex}\\
(b) Every $\We$-orbit of $\RkX$ contains exactly one special face of $X$.
\end{cor}

By part (b) of the last corollary, every face $R$ can be parametrized in the form $R=\sigma R(\Th)$ with $\sigma\in\We$, and with a uniquely determined 
special set $\Th$, which we call the {\it type} of $R$. The element $\sigma$ is uniquely determined, if we restrict to the minimal coset 
representatives $\We^{\Th\cup\Th^\bot}$ of $\We/\We_{\Th\cup\Th^\bot}$. The next theorem shows how to express the inclusion and the lattice operations by this 
parametrization. The equivalence of part (a) has already been given as Proposition 2.5 of \cite{M1}. Part (b) is new. 

To state this theorem define a function $red:\We\to I$ in the following way: Set $\red{1}:=\emptyset$. If $\sigma\in\We\setminus\{1\}$ and $\sigma_{i_1}\cdots \sigma_{i_k}$ is a reduced expression for $\sigma$ with simple reflections $\sigma_{i_1},\,\ldots,\,\sigma_{i_k}$ corresponding to $i_1,\,\ldots,\,i_k\in I$ set
\begin{eqnarray*}
\red{\sigma}:= \{i_1,\,\ldots,\,i_k\}\subseteq I.
\end{eqnarray*} 
The function $red$ is well defined, because this subset of $I$ is independent of the chosen reduced expression by a result of J. Tits, compare for example Theorem 2.11 of \cite{Ro}.
Denote by $\leq$ the Bruhat order on $\We$. The function $red$ is order preserving, i.e.,  $\sigma\leq \ti{\sigma}$ implies 
$\red{\sigma}\subseteq \red{\ti{\sigma}}$, $\sigma,\,\ti{\sigma}\in\We$. 

Recall the following convention to cut short many formulas: We often identify in the obvious way the set of indices $I$ with the set of simple roots. 

For the theorem note also that the lattice intersection of two faces coincides with the set theoretical intersection of the faces.

\begin{thm}\label{FL4}\mb{}
(a) Let $\sigma,\sigma'\in\We$, $\Th,\Th'$ be special. Then:
\begin{eqnarray}\label{sub}
   \sigma'R(\Th') \subseteq \sigma R(\Th)\quad \iff \quad \Th' \supseteq \Th \mb{ and } \sigma^{-1}\sigma' \in \We_{\Th^\bot}\We_{\Th'}.  
\end{eqnarray}
Different faces of $X$ of the same type are not comparable by $\subseteq$.\vspace*{1ex}\\
(b) Let $\tau\in\mb{}^{\Th_1\cup\Th_1^\bot}\We^{\Th_2\cup\Th_2^\bot}$, $\Th_1,\Th_2$ be special. Then
\begin{eqnarray}
   R(\Th_1)\cap \tau R(\Th_2)   &=&  R(\Th_1\cup\Th_2\cup \red{\tau}), \label{int}\\
   R(\Th_1)\vee\tau R(\Th_2)   &=&  R((\Th_1\cap \tau\Th_2)^\infty). \label{join}
\end{eqnarray}
\end{thm}
\begin{rem}\label{intjoinpar} The lattice intersection and lattice join of two arbitrary faces can be reduced easily to the formulas (\ref{int}), (\ref{join}). Use  
the formulas for the stabilizers of the special faces. 

In particular, the following slightly more general formulas than (\ref{int}) and (\ref{join}) hold:
Let $\tau\in\mb{}^{\Th_1^\bot}\We^{\Th_2^\bot}$, $\Th_1,\Th_2$ be special. Then
\begin{eqnarray*}
   R(\Th_1)\cap \tau R(\Th_2)   =  R(\Th_1\cup\Th_2\cup \red{\tau}). 
\end{eqnarray*}
Let $\tau\in\mb{}^{\Th_1}\We^{\Th_2}$, $\Th_1,\Th_2$ be special. Then
\begin{eqnarray*}
   R(\Th_1)\vee\tau R(\Th_2)   =  R((\Th_1\cap \tau\Th_2)^\infty).
\end{eqnarray*}
\end{rem}
\Proof The equivalence of part (a) has been given in Proposition 2.5 of \cite{M1}. From this equivalence and formula (\ref{Stwhole}) follows immediately that different faces of 
the same type are not comparable. 

We now prove formula (\ref{int}) of part (b):
We first show that the set
\begin{eqnarray}\label{DefTh}
         \Th:=\Th_1\cup\Th_2\cup\red{\tau}
\end{eqnarray} 
is a special subset of $I$. Obviously, $\tau\in\We_\Th$ by the definition of $\Th$. Decompose $\tau$ in the form $\tau=\tau_\infty\tau_0$ with $\tau_\infty\in\We_{\Th^\infty}$, $\tau_0\in\We_{\Th^0}$.
The sets $\Th^\infty$ and $\Th^0$ consist of connected components of $\Th$. Thus $\Th^0\subseteq(\Th^\infty)^\bot$.
Furthermore, $\Th_2\subseteq\Th$ implies $\Th_2\subseteq \Th^\infty$, which in turn implies $(\Th_2)^\bot \supseteq (\Th^\infty)^\bot$. Therefore,
\begin{eqnarray*}
   \We_{\Th_1\cup\Th_1^\bot} \tau \We_{\Th_2\cup\Th_2^\bot}   =   \We_{\Th_1\cup\Th_1^\bot} \tau_\infty \tau_0 \We_{\Th_2\cup\Th_2^\bot}    
                                                              =   \We_{\Th_1\cup\Th_1^\bot} \tau_\infty \We_{\Th_2\cup\Th_2^\bot}.
\end{eqnarray*}
Since $\tau$ is an element of minimal length in this double coset we get $l(\tau_\infty)\geq l(\tau)=l(\tau_\infty)+l(\tau_0)$. This implies $l(\tau_0)=0$. Therefore,
$\tau_0=1$ and $\tau=\tau_\infty$. Because of $\Th_1\subseteq \Th^\infty$, $\Th_2\subseteq \Th^\infty$, and $\red{\tau}\subseteq\Th^\infty$ we get
\begin{eqnarray*}
    \Th^\infty\subseteq \Th=\Th_1\cup\Th_2\cup\red{\tau}\subseteq \Th_1\cup\Th_2\cup\Th^\infty =\Th^\infty,
\end{eqnarray*} 
which shows that $\Th$ is special.

The intersection of the faces $R(\Th_1)$ and $\tau R(\Th_2)$ is again a face of $X$, i.e., $R(\Th_1)\cap \tau R(\Th_2)= \eta R(J)$ for some $\eta\in\We$ and some $J\subseteq I$ 
special. By part (a) we find
\begin{eqnarray*}
   J\supseteq \Th_1\cup\Th_2  \;\mb{ and }\;\eta \in \We_{\Th_1^\bot}\We_J,\; \tau^{-1}\eta\in \We_{\Th_2^\bot}\We_J.
\end{eqnarray*}
It follows $\tau\in \We_{\Th_1^\bot} \We_J\We_{\Th_2^\bot}$. Therefore, there exists an element 
\begin{eqnarray*}
   x\in \We_J \cap \We_{\Th_1^\bot} \tau \We_{\Th_2^\bot}\subseteq  \We_{\Th_1\cup \Th_1^\bot} \tau \We_{\Th_2\cup \Th_2^\bot}.
\end{eqnarray*}
Since $\tau$ is a minimal double coset representative we have $\tau\leq x $. Therefore $\red{\tau}\subseteq \red{x} \subseteq J$. Because of $\Th_1\cup\Th_2\subseteq J$ it follows $\Th\subseteq J$.

Only by definition (\ref{DefTh}) and part (a) we get $R(\Th)\subseteq R(\Th_1)\cap \tau R(\Th_2)=\eta R(J)$ immediately. Applying (a) once more shows
the reverse inclusion $\Th\supseteq J$. 

It follows $R(\Th)\subseteq \eta R(\Th)=R(\Th_1)\cap \tau R(\Th_2)$, which implies $R(\Th)=\eta R(\Th)=R(\Th_1)\cap \tau R(\Th_2)$ by part (a).

Finally, we prove formula (\ref{join}) of part (b): We first show the inclusion
\begin{eqnarray}\label{veeincl}
  R(\Th_1)\vee\tau R(\Th_2) \subseteq R((\Th_1\cap \tau\Th_2)^\infty) .
\end{eqnarray}
Obviously, $(\Th_1\cap \tau\Th_2)^\infty\subseteq \Th_1$. It follows $R((\Th_1\cap \tau\Th_2)^\infty)\supseteq R(\Th_1)$. Now let $\la\in \tau R(\Th_2)$. Then it holds 
\begin{eqnarray*}
   0=(\tau^{-1}\la)(h_i)=\la(\tau h_i)=\la(h_{\tau \al_i})\;\mb{ for all } i\in\Th_2.
\end{eqnarray*}
In particular, $\la(h_i)=0$ for all $i\in \Th_1\cap\tau \Th_2$. Therefore, we get $\tau R(\Th_2)\subseteq R((\Th_1\cap \tau\Th_2)^\infty)$.

The join of the faces $R(\Th_1)$ and $\tau R(\Th_2)$ is a face, i.e., $R(\Th_1)\vee \tau R(\Th_2)=\eta R(J)$ for some $\eta\in\We$ and some $J\subseteq I$ 
special. By part (a) we get
\begin{eqnarray*}
     J\subseteq \Th_1\cap\Th_2  \;\mb{ and }\;\eta \in \We_{\Th_1}\We_{J^\bot},\; \tau^{-1}\eta\in \We_{\Th_2}\We_{J^\bot}.
\end{eqnarray*}
It follows $\tau\in \We_{\Th_1} \We_{J^\bot}\We_{\Th_2}$. Therefore, there exists an element 
\begin{eqnarray*}
   x\in \We_{J^\bot} \cap \We_{\Th_1} \tau \We_{\Th_2}\subseteq  \We_{J^\bot} \cap  \We_{\Th_1\cup \Th_1^\bot} \tau \We_{\Th_2\cup \Th_2^\bot}.
\end{eqnarray*}
Since $\tau$ is a minimal double coset representative we have $\tau\leq x $. This implies $\red{\tau}\subseteq \red{x} \subseteq J^\bot$. Therefore, $\tau J=J$. It follows $J\subseteq \Th_1\cap \tau\Th_2$. Since $J$ is special we get $J\subseteq (\Th_1\cap \tau\Th_2)^\infty$, from which follows by part (a)
\begin{eqnarray*}
  \eta R( (\Th_1\cap \tau\Th_2)^\infty)\subseteq \eta R(J)=R(\Th_1)\vee \tau R(\Th_2).
\end{eqnarray*}
If we use the inclusion (\ref{veeincl}) and part (a) once more we get $R(\Th_1)\cap \tau R(\Th_2)=R((\Th_1\cap \tau\Th_2)^\infty)$.
\qed 

Later, we also need the following proposition, stated as Proposition 2.7 of \cite{M1}. It shows how the face lattice of the Tits cone $X_{A_J}$ associated to a generalized Cartan submatrix $A_J$ of the generalized Cartan matrix $A$ can be found as a sublattice the face lattice of the Tits cone $X$ associated to $A$.
\begin{prop}\label{FL5}
Let $\emptyset\neq J\subseteq I$. Identify the Coxeter group $\We(A_J)$ with the parabolic subgroup 
$\We_J$ of $\We$. We get an embedding of lattices 
\begin{eqnarray*}
  \Upsilon_J:\;\; Fa(X_{A_J})\, \to      \RkX     
\end{eqnarray*}
by $\Upsilon_J(\tau R_{A_J}(\Th)):= \tau R(\Th)$, where $\tau\in \We(A_J)=\We_J$, and $\Th\subseteq J$ special. Its image is given by 
\begin{eqnarray*}
  \RkX_J &:=&  \Mklz{R\in\RkX}{R\supseteq R(J^\infty)}\\ 
         &=& \Mklz{\tau R(\Th)}{\tau\in\We_J,\;\Th\subseteq J \mb{ special }}\\
         &=&\Mklz{\tau R(\Th)}{\tau\in \We_{J^\infty},\;\Th\subseteq J^\infty \mb{ special }}.
\end{eqnarray*}
\end{prop}
To simplify the notation of many formulas we set $\RkX_\emptyset:=\{X\}$.\vspace*{1ex}

If $R$ is a face of $X$ and $\overline{F}$ is a closed facet of $X$, then $R\cap \overline{F}$ is a face of $\overline{F}$. Therefore, $R\cap \overline{F}$ is also a closed facet of $X$. Explicitely, we have the following Theorem:
\begin{thm}\label{FL6} Let $\Th$ be special subset of $I$, and $J\subseteq I$. Let $\sigma\in\mb{}^{\Th\cup\Th^\bot}\We^J$. Then
\begin{eqnarray*}
  R(\Th)\cap\sigma \overline{F}_J= \overline{F}_{\Th\cup J\cup\,\red{\sigma}}.
\end{eqnarray*}
\end{thm}
\begin{rem} The formula for the intersection of an arbitrary face with an arbitrary facet can be derived easily form this formula. Use the formulas for the stabilizers of the special faces and of the closed facets as a whole. 
\end{rem}
\Proof Note first that $\overline{F}_{\Th\cup J\cup\,\red{\sigma}}=\sigma \overline{F}_{\Th\cup J\cup\,\red{\sigma}}$. 
Since $R(\Th)=\We_{\Th^\bot}\overline{F}_\Th\supseteq \overline{F}_{\Th\cup J\cup\,\red{\sigma}}$ and $\overline{F}_J\supseteq \overline{F}_{\Th\cup J\cup\,\red{\sigma}}$ the inclusion ``$\supseteq$'' of the Theorem follows immediately. We now show the reverse inclusion ``$\subseteq$''. The intersection of $R(\Th)$ and $\sigma\overline{F}_J$ is a face of $\sigma\overline{F_J}$, i.e., 
\begin{eqnarray*} 
    R(\Th)\cap \sigma\overline{F}_J=\sigma\overline{F}_L
\end{eqnarray*} 
for some $L\supseteq J$. The inclusion $\sigma F_L\subseteq R(\Th)=\We_{\Th^\bot}\overline{F}_\Th$ implies $\Th\subseteq L$ and 
$\sigma F_L\subseteq \We_{\Th^\bot} F_L$,
from which in turn follows 
\begin{eqnarray}\label{FR1}
\sigma\in\We_{\Th^\bot}\We_L= \We_{\Th^\bot}\mb{}^{\Th^\bot\cap L}(\We_L)   .
\end{eqnarray} 
Furthermore, 
\begin{eqnarray}\label{FR2}
 \mb{}^{\Th^\bot\cap L}(\We_L)\subseteq \mb{}^{\Th^\bot}\We,
\end{eqnarray}
which can be seen as follows. Let $\tau\in \mb{}^{\Th^\bot\cap L}(\We_L)$. If $i\in \Th^\bot\cap L$ then trivially $\tau \al_i\in \prW$. If $i\in\Th^\bot\setminus L$ then
\begin{eqnarray*}
\tau \al_i=\al_i+ \mb{ a linear combination in the elements }\al_l,\; l\in L, 
\end{eqnarray*}
from which follows $\tau \al_i\in \prW$.

Because of $\sigma\in \mb{}^{\Th\cup\Th^\bot}\We^J\subseteq \mb{}^{\Th\cup\Th^\bot}\We\subseteq \mb{}^{\Th^\bot}\We$, from (\ref{FR1}) and (\ref{FR2}) follows $\sigma\in \mb{}^{\Th^\bot\cap L}(\We_L)$. In particular, $\red{\sigma}\subseteq L$.

Collecting the inclusions $J\subseteq L$, $\Th\subseteq L$, and $\red{\sigma}\subseteq L$ we find $\Th\cup J\cup\red{\sigma}\subseteq L$, from which follows
\begin{eqnarray*}
  R(\Th)\cap \sigma\overline{F}_J=\sigma\overline{F}_L\subseteq \sigma \overline{F}_{\Th\cup J\cup\,\red{\sigma}}=\overline{F}_{\Th\cup J\cup\,\red{\sigma}}.
\end{eqnarray*}
\qed  
%
%
%%%%%%%%%%%%%%%%%%%%%%%%%%%%%%%%%%%%%%%%%%%%%%%%%%%%%%%%%%%%%%%%%%%%%%%%%%%%%%%%%%%%%%%%%%%%%%%%%%%%%%%%%%%%%%%%%%%%%%%%%%%%%%%%%%%%%%%%
%
%
\section{The face monoid $\GD$ associated to a Kac-Moody group $G$, and the face monoid $\WeD$ associated to its Weyl group $\We$}
The face monoid $\GD$ of \cite{M1} is a very natural infinite-dimensional algebraic monoid. It has similar properties as a reductive algebraic monoid. Its unit group $G$ coincides, up to a slightly extended maximal torus, with the special Kac-Moody group defined representation theoretically in \cite{KP1}. Therefore, we call $\GD$ the face monoid associated to the Kac-Moody group $G$. This monoid is a purely non-classical phenomenon. In the classical case, i.e., $G$ a semisimple simply connected algebraic group, $\GD$ coincides with $G$ itself. In this section we review the algebraic properties of $\GD$ which are used later. As new results we obtain an elegant explicit description for the multiplication of two elements of $\GD$ and a version of a monoid theoretical construction of $\GD$.
%%%%%%%%%%%%%%%%%%%%%%%%%%%%%%%%%%%%%%%%%%%%%%%%%%%%%%%%%%%%%%%%%%%%%%%%%%%%%%%%%%%%%%%%%%%%%%%%%%%%%%%%%%%%%%%%%%%%%%%%%%%%%%%%%%%
%Please add some additional space:
\vspace*{1ex}
%
%%%%%%%%%%%%%%%%%%%%%%%%%%%%%%%%%%%%%%%%%%%%%%%%%%%%%%%%%%%%%%%%%%%%%%%%%%%%%%%%%%%%%%%%%%%%%%%%%%%%%%%%%%%%%%%%%%%%%%%%%%%%%%%%%%%%
%

The face monoid $\GD$ can be defined quickly as follows: Recall that we denote by $P$ a weight lattice, by $P^+$ the dominant weights of $P$. Recall that $L(\La)$ denotes an irreducible highest weight module of highest weight $\La$, and $P(\La)$ its set of weights. Now the Kac-Moody group $G$ acts faithfully on $\bigoplus_{\La\in P^+}L(\La)$. Identify $G$ with the corresponding subgroup of the monoid $End(\bigoplus_{\La\in P^+}L(\La))$ of linear endomorphisms of $\bigoplus_{\La\in P^+}L(\La)$. For every face $R$ of the Tits cone $X$ we get a projection $e(R)\in End(\bigoplus_{\La\in P^+}L(\La))$ by
\begin{eqnarray}\label{defeR}
  e(R) v_\la:=\left\{\begin{array}{cl}
  v_\la   & \mb{ if } \la\in R\\
    0     & \mb{ else}
                        \end{array}\right. ,\; \mb{ where }\;v_\la\in L(\La)_\la,\;\la\in P(\La),\;\La\in P^+.
\end{eqnarray} 
The face monoid $\GD$ is the submonoid of $End(\bigoplus_{\La\in P^+}L(\La))$ generated by the Kac-Moody group $G$ and the elements $e(R)$, $R$ a face of the Tits cone $X$.

This is a quick definition of $\GD$. It looks a bit strange because $\GD$ is defined in an ad hoc way as a submonoid of this very big monoid $End(\bigoplus_{\La\in P^+}L(\La))$. In \cite{M1} we obtained $\GD$ by a Tannaka reconstruction. Roughly, this shows that $\GD$ is the smallest monoid acting "reasonably" on the irreducible highest weight modules $L(\La)$, $\La\in P^+$, compatible with their integrable duals.
%
%%%%%%%%%%%%%%%%%%%%%%%%%%%%%%%%%%%%%%%%%%%%%%%%%%%%%%%%%%%%%%%%%%%%%%%%%%%%%%%%%%%%%%%%%%%%%%%%%%%%%%%%%%%%%%%%%%%%%%%%%%%%%%%%%%%
%Please add some additional space:
\vspace*{1ex}
%
%%%%%%%%%%%%%%%%%%%%%%%%%%%%%%%%%%%%%%%%%%%%%%%%%%%%%%%%%%%%%%%%%%%%%%%%%%%%%%%%%%%%%%%%%%%%%%%%%%%%%%%%%%%%%%%%%%%%%%%%%%%%%%%%%%%%
%

The following theorem gives the $G\times G$-orbit decomposition of $\GD$. In addition, it gives an elegant algebraic description of the monoid $\GD$. Part (a) and (b) have been given in \cite{M1} as Proposition 2.31 and Theorem 2.20 (a). The explicit formula for the multiplication of two elements of $\GD$ in part (c) is new. 

To state the theorem we first introduce some notation. Let $\Th$ be special. Because of the Levi decomposition 
\begin{eqnarray*}
P_{\Th\cup\Th^\bot}= L_{\Th\cup\Th^\bot}\ltimes U^{\Th\cup\Th^\bot},
\end{eqnarray*}
and the decomposition 
\begin{eqnarray*}
 L_{\Th\cup\Th^\bot}=(G_\Th\times G_{\Th^\bot})\rtimes T^{\Th\cup\Th^\bot}, \quad\mb{ where }\quad 
    T^{\Th\cup\Th^\bot}:= T_{I\setminus (\Th\cup\Th^\bot)}T_{rest}, 
\end{eqnarray*}
we get the decompositions
\begin{eqnarray*}
    P_{\Th\cup\Th^\bot}  &= &\left(G_{\Th^\bot}\rtimes T^{\Th\cup\Th^\bot}\right)\ltimes \left(G_\Th\ltimes U^{\Th\cup\Th^\bot}\right),\\
   \left(P_{\Th\cup\Th^\bot}\right)^- &=& \left((U^{\Th\cup\Th^\bot})^-\rtimes G_\Th\right) \rtimes\left(G_{\Th^\bot}\rtimes T^{\Th\cup\Th^\bot}\right).
\end{eqnarray*}
The projections belonging to these semidirect products are denoted by
\begin{eqnarray*}
 p_\Th:=p_\Th^+\,&:&\; P_{\Th\cup\Th^\bot}  \to  G_{\Th^\bot}\rtimes T^{\Th\cup\Th^\bot},\\
 p^-_\Th\,&:&\; \left(P_{\Th\cup\Th^\bot}\right)^- \to  G_{\Th^\bot}\rtimes T^{\Th\cup\Th^\bot}.
\end{eqnarray*}

Denote by $\widetilde{\quad}:\We\to N$ the cross section of the canonical map from $N$ to $\We$ defined by $\widetilde{1}=1$, $\widetilde{\sigma_i}=n_i$, $i\in I$, and
$\widetilde{\sigma\tau}=\widetilde{\sigma}\widetilde{\tau}$ if $l(\sigma\tau)=l(\sigma)+l(\tau)$, $\sigma,\tau\in\We$.
\begin{thm}\label{FM1} (a) We have
\begin{eqnarray*}
  \GD = \dot{\bigcup_{\Th\;special}} G e(R(\Th)) G  .
\end{eqnarray*}
(b) Let $g_1,g_2, h_1, h_2\in G$, and let $\Th_1$, $\Th_2$ be special. The following statements are equivalent:
\begin{itemize}
\item[(i)] $g_1 e(R(\Th_1)) h_1= g_2 e(R(\Th_2)) h_2 $.
\item[(ii)] $\Th_2 =\Th_1$ and $(g_2)^{-1}g_1\in P_{\Th_1\cup\Th_1^\bot}$, $h_2(h_1)^{-1}\in\left(P_{\Th_1\cup\Th_1^\bot}\right)^-$ with $p_{\Th_1}((g_2)^{-1}g_1) = p_{\Th_1}^-(h_2 (h_1)^{-1})$.
\end{itemize}
(c) Let $g_1,g_2, h_1, h_2\in G$, and let $\Th_1$, $\Th_2$ be special. Consider the elements $g_1 e(R(\Th_1)) h_1$, $g_2 e(R(\Th_2)) h_2 $. By the generalized Birkhoff decomposition write $h_1g_2$ in the form
\begin{eqnarray*}
  h_1 g_2= a_1 \widetilde{\sigma} a_2
\end{eqnarray*}
where $a_1\in P_{\Th_1\cup\Th_1^\bot}^-$, $a_2\in  P_{\Th_2\cup\Th_2^\bot}$, and $\sigma\in \mb{}^{\Th_1\cup\Th_1^\bot}\We^{\Th_2\cup\Th_2^\bot}$. Then
\begin{eqnarray*}
 g_1 e(R(\Th_1)) h_1 \, g_2 e(R(\Th_2)) h_2 = g_1 p_{\Th_1}^-(a_1)\, e(\,R(\Th_1\cup\Th_2\cup\red{\sigma})\,)\,p_{\Th_2}(a_2)h_2.
\end{eqnarray*}
\end{thm}
\Proof It remains to show (c). By applying two times part (b) we find
\begin{eqnarray*}
      && g_1 e(R(\Th_1)) h_1 g_2 e(R(\Th_2)) h_2 =  g_1 e(R(\Th_1)) a_1 \widetilde{\sigma} a_2 e(R(\Th_2)) h_2 \\
      && = g_1 p_{\Th_1}^-(a_1) e(R(\Th_1)) \widetilde{\sigma} e(R(\Th_2)) p_{\Th_2}(a_2)h_2.
\end{eqnarray*}
Now consider the middle part $e(R(\Th_1)) \widetilde{\sigma} e(R(\Th_2))$ of this expression. By Theorem 2.20 (b) of \cite{M1} and by formula (\ref{int}) for the 
intersection of two faces in Theorem \ref{FL4} we get
\begin{eqnarray*}
    e(R(\Th_1)) \widetilde{\sigma} e(R(\Th_2)) =  e(R(\Th_1)\cap \sigma R(\Th_2))\widetilde{\sigma}=e(R(\Th_1\cup\Th_2\cup\red{\sigma}))\widetilde{\sigma}.
\end{eqnarray*}
By (b) the last expression coincides with $e(R(\Th_1\cup\Th_2\cup\red{\sigma}))$.\qed 

There are also Bruhat and Birkhoff decompositions for $\GD$, which we want to state next. Denote by $\ND$ the monoid generated by $N$ and the elements 
$e(R)$, $R$ a face of the Tits cone $X$. We get a congruence relation on $\ND$ as follows:
\begin{eqnarray*}
\quad\widehat{n}\:\sim\:\widehat{n}'\quad:\iff\quad\widehat{n}T\:=\:\widehat{n}'T\quad\iff \quad\widehat{n}'\:\in\:\widehat{n}T\quad\iff\quad\widehat{n}\in\widehat{n}'T
\end{eqnarray*}
Denote the monoid $\ND/\sim$ by $\ND/T$. Note that the monoid $\ND/T$ contains the Weyl group $N/T\cong \We$.

The following Theorem has been given as Theorem 2.33 in \cite{M1}.
\begin{thm}\label{FM2} There are the Bruhat and Birkhoff decompositions 
\begin{eqnarray*}
   \GD = \dot{\bigcup_{\widehat{w}\,\in\,\widehat{N}/T}} B^\eps\widehat{w}B^\delta =
  \dot{\bigcup_{\widehat{n}\,\in\,\widehat{N}}} U^\eps\widehat{n}U^\delta \; ,\qquad\qquad\qquad\qquad\qquad\qquad\qquad\\
        \qquad\qquad\qquad \mb{ where } \quad  (\epsilon,\delta)\;=\;\underbrace{(+,+)\;,\;(-,-)}_{Bruhat}\;,\;\underbrace{(+,-)\;,\;(-,+)}_{Birkhoff}\;\;.
\end{eqnarray*}
\end{thm}

We now want to describe the monoid $\ND/T$. We do this in two steps. We first define and describe a monoid $\WeD$ abstractly, which then identifies with $\ND/T$.\vspace*{1ex}

The Weyl group $\We$ acts on the monoid (\,$\RkX\,,\,\cap$\,). The semidirect product $\We\ltimes \RkX $ consists of the set 
$\We\times \RkX$ with the structure of a monoid given by 
\begin{eqnarray*}
  (\sigma, R)\cdot (\tau, S) := (\sigma \tau, \tau^{-1} R\cap S) .
\end{eqnarray*}
It is easy to see that we get a congruence relation on $\We\ltimes \RkX $ by
\begin{eqnarray*}
   (\sigma, R) \sim (\sigma', R') &:\iff & R\:=\:R'\quad\mb{and}\quad \sigma^{-1}\sigma'\in Z_{\mathcal W}(R)\;.
\end{eqnarray*}
We denote the congruence class of $(\sigma, R)$ by $\sigma \ve{R}$.
\begin{defn} We call the monoid $\widehat{\We} := (\We\ltimes \RkX )/\sim $ the face monoid associated to the Weyl group $\We$. 
(In \cite{M1} it was called Weyl monoid.)  
\end{defn}

It is easy to check, that we get embeddings of monoids
\begin{eqnarray*}
  \begin{array}{ccc}      
    \RkX &\to     &\WeD\\
      R  &\mapsto &1 \ve{R}
  \end{array} 
  & \quad , \quad &
  \begin{array}{ccc}
         \We &\to     & \WeD   \\                       
      \sigma &\mapsto & \sigma \ve{X}
  \end{array}\;. 
\end{eqnarray*}
We denote the element $1 \ve{R}$ by $\ve{R}$, and the image of $\RkX$ by ${\mathcal E}$. We denote the element $\sigma \ve{X}$ by $\sigma$, i.e., we identify the Weyl group $\We$ with its image in $\WeD$.

The following two Propositions have been given as Proposition 2.24 and 2.25 of \cite{M1}.
\begin{prop}\label{FM3} $\We$ is the unit group of $\WeD$, ${\mathcal E}$ is the set of idempotents of $\WeD$, and 
$\WeD$ is a unit regular, i.e., $\WeD=\We \,{\mathcal E}= {\mathcal E}\,\We$.
\end{prop}
\begin{prop}\label{FM3b} $\WeD$ is an inverse monoid. The inverse map $\mb{}^{inv}:\WeD\to\WeD$ is given by 
\begin{eqnarray*}
 (\sigma \ve{R})^{inv}=\ve{R}\sigma^{-1} \quad \mb{ where }\quad  \sigma\in\We,\; R \mb{ a face of }X.
\end{eqnarray*}
\end{prop}
The group $\We\times \We$ acts in the obvious way on $\WeD$. It is trivial to show: 
\begin{prop} The $\We\times \We$-orbit decomposition of $\WeD$ is 
\begin{eqnarray*}
\WeD=\dot{\bigcup_{\Th\,special}} \We\ve{R(\Th)}\We.
\end{eqnarray*}
\end{prop}
Now the following Proposition, which has been given as Proposition 2.30 of \cite{M1}, relates the monoid $\ND/T$ and the monoid $\WeD$.
\begin{prop}\label{FM4}
The monoid $\ND/T$ is isomorphic to $\WeD$ by the following map:
\begin{eqnarray*}
  \kappa\,:\quad\WeD\;\:\;    & \to     &   \;\;\ND/T       \\
       \sigma \ve{R}        & \mapsto & \widetilde{\sigma} e(R)T
\end{eqnarray*}
\end{prop}
We identify $\ND/T$ with $\WeD$ by this map.\vspace*{1ex}

Next we want to mention some normal forms of the elements of $\WeD$ and $\GD$ adapted to the $\We\times\We$-orbit decomposition of $\WeD$ and the Bruhat decompositions of $\GD$. The following proposition, which is very easy to check, has been given as a part of Proposition 6 in \cite{M3}.
\begin{prop}\label{FM5} Let $\Th$ be special. Let $\hat{\sigma}\in \We\ve{R(\Th)}\We$. It holds:
\begin{itemize}\item[(I)] The element $\hat{\sigma}$ can be uniquely written in the form
\begin{eqnarray*}
 \hat{\sigma}=\sigma_1\ve{R(\Th)}\sigma_2 \quad\mb{ with }\quad\sigma_1\in\We^\Th,\;\;\sigma_2\in\mb{}^{\Th\cup \Th^\bot}\We.
\end{eqnarray*}
\item[(II)] The element $\hat{\sigma}$ can be uniquely written in the form
\begin{eqnarray*}
 \hat{\sigma}=\sigma_1\ve{R(\Th)}\sigma_2 \quad\mb{ with }\quad\sigma_1\in\We^{\Th\cup\Th^\bot},\;\;\sigma_2\in\mb{}^\Th\We.
\end{eqnarray*}
\end{itemize}
\end{prop}
We call these forms {\it normal form I} and {\it normal form II}.\vspace*{1ex}

The following theorem is the underlying algebraic part for the description of the Bruhat cells of $\GD$ given in Theorem 14, Part (1), of \cite{M3}.
\begin{thm}\label{FM6} Let $\hat{n}\in\ND$ and $\hat{g}\in U\hat{n}U$. Let $\hat{\sigma}\in\WeD=\ND/T$ be the projection of $\hat{n}$. 
Let $\hat{\sigma}=\sigma_1\ve{R(\Th)}\sigma_2$ be the normal form I of $\hat{\sigma}$. The element $\hat{g}$ can be uniquely written in the form
\begin{eqnarray*}
 \hat{g}=u_1\hat{n} u_2 \quad\mb{ with }\quad u_1\in U\cap \sigma_1 (U^\Th)^- (\sigma_1)^{-1},\;\;u_2\in U\cap (\sigma_2)^{-1}U^\Th\sigma_2.
\end{eqnarray*}
Furthermore, $U_{\sigma_1}:=U\cap \sigma_1 U^- (\sigma_1)^{-1}=U\cap \sigma_1 (U^\Th)^- (\sigma_1)^{-1}$.\vspace*{1ex}

Let $\hat{n}\in\ND$ and $\hat{g}\in U^-\hat{n}U^-$. Let $\hat{\sigma}\in\WeD=\ND/T$ be the projection of $\hat{n}$. 
Let $\hat{\sigma}=\sigma_1\ve{R(\Th)}\sigma_2$ be the normal form II of $\hat{\sigma}$. The element $\hat{g}$ can be uniquely written in the form
\begin{eqnarray*}
 \hat{g}=u_1\hat{n} u_2 \quad\mb{ with }\quad u_1\in U ^-\cap \sigma_1 (U^\Th)^- (\sigma_1)^{-1},\;\;u_2\in U^-\cap (\sigma_2)^{-1}U^\Th\sigma_2.
\end{eqnarray*}
Furthermore, $ U_{(\sigma_2)^{-1}}^-:=U^-\cap (\sigma_2)^{-1}U\sigma_2=U^-\cap (\sigma_2)^{-1}U^\Th\sigma_2$.

\end{thm}
The monoid $\GD$ can be obtained abstractly from the Kac-Moody group $G$ and the monoid $(Fa(X),\cap)$. There are three related versions which use: 
\begin{itemize}
\item[(i)] The root groups data system $(U_\al)_{\al\in\Delta_{re}}$, $T$.
\item[(ii)] Some subgroups of $G$.
\item[(iii)] The projections $p_{\Th}$, $p_{\Th}^*$, $\Th$ special.
\end{itemize}
In  Theorem 2.34 of \cite{M1} only version (i) has been given. You obtain version (ii) resp. (iii) easily. Replace 
the formulas (2) of Theorem 2.34 of \cite{M1} by the formulas ($\al$), ($\beta$), ($\gamma$) resp. ($\delta$) of the proof of this theorem.

Now in Theorem \ref{FL4} (b) the intersection of the faces of the Tits cone $X$ has been reduced to the intersection of the special faces of $X$. This can be used to rewrite Theorem 2.34 of \cite{M1} to obtain the monoid $\GD$ abstractly from the Kac-Moody group $G$ and the monoid $(\{\mb{ Special sets of $I$ }\},\,\cup)$, which is isomorphic to the submonoid 
\begin{eqnarray*}
    E_{sp}:=\{\,e(R(\Th))\,\mid\,\Th \mb{ a special set of }I\,\}
\end{eqnarray*} 
of $\GD$:
\begin{thm} Let $\sim$ be the congruence relation on the free product $G*E_{sp}$ of the monoids $G$ and $E_{sp}$ generated by the following relations:
\begin{itemize}
\item[(1)] For every special set $\Th\subseteq I$ and every $t\in T$:  $\quad te(R(\Th))t^{-1}\sim e(R(\Th))$.
\item[(2)] For every special set $\Th\subseteq I$, every $x\in U_\al, \al\in \rW$:
\begin{itemize}
\item[(a)] If $\al\in \We_{\Th\cup\Th^\bot}(\Th\cup\Th^\bot)$: $\quad xe(R(\Th))x^{-1}\sim e(R(\Th))$.
\item[(b)] If $\al\in \prW\setminus(\We_{\Th^\bot}\Th^\bot)\cup\We_{\Th}\Th$: $\quad xe(R(\Th))\sim e(R(\Th))$.
\item[(c)] If $\al\in \nrW\setminus(\We_{\Th^\bot}\Th^\bot)\cup\We_{\Th}\Th$: $\quad e(R(\Th))x\sim e(R(\Th))$.
\end{itemize}
\item[(3)] For every special set $\Th\subseteq I$, $\Xi\subseteq I$, for every $\sigma\in\mb{}^{\Th\cup\Th^\bot}\We^{\Xi\cup\Xi^\bot}$, for every $n_\sigma\in N$ projecting to $\sigma$:  $\quad e(R(\Th))n_\sigma e(R(\Xi))n_\sigma^{-1}\sim e(R(\Th\cup\Xi\cup \red{\sigma}))$.
\end{itemize}
\end{thm}
Then $\GD$ is isomorphic in the obvious way to $(G*E_{sp})/\sim$.
\begin{rem} Because of (2) (a), the congruence in (3) may be replaced by 
\begin{eqnarray*}
    n_\sigma^{-1}e(R(\Th))n_\sigma e(R(\Xi))\sim e(R(\Th\cup\Xi\cup \red{\sigma})).
\end{eqnarray*} 
Because of (1) and  (2) (c), it may also be replaced by 
\begin{eqnarray*}
    e(R(\Th))\widetilde{\sigma}(R(\Xi))\sim e(R(\Th\cup\Xi\cup \red{\sigma})).
\end{eqnarray*}
\end{rem}
\section{Actions of the face monoid $\WeD$ associated to $\We$ on the Coxeter complex $\mathcal C$}
The Weyl group $\We$ acts naturally on its Coxeter complex $\mathcal C$. We first describe an extension to an action of the face monoid $\WeD$, obtained  directly by representation theory. Unfortunately, this action does not fit to the structure of the Coxeter complex $\mathcal C$, because it does not preserve its natural order. We therefore call this action the bad action. At this point one could think that the face monoid $\WeD$ does not fit properly to the Coxeter complex $\mathcal C$. But this is not the case: We show how to modify the bad action to obtain two different order preserving actions of $\WeD$ on the Coxeter complex $\mathcal C$. Unexpectingly, there exists more than one such action. 
%
%%%%%%%%%%%%%%%%%%%%%%%%%%%%%%%%%%%%%%%%%%%%%%%%%%%%%%%%%%%%%%%%%%%%%%%%%%%%%%%%%%%%%%%%%%%%%%%%%%%%%%%%%%%%%%%%%%%%%%%%%%%%%%%%%%%
%Please add some additional space:
\vspace*{1ex}
%
%%%%%%%%%%%%%%%%%%%%%%%%%%%%%%%%%%%%%%%%%%%%%%%%%%%%%%%%%%%%%%%%%%%%%%%%%%%%%%%%%%%%%%%%%%%%%%%%%%%%%%%%%%%%%%%%%%%%%%%%%%%%%%%%%%%%
%

To describe the stabilizers of our actions we need the analogues of the parabolic subgroups of $\We$, the parabolic submonoids of $\WeD$, which we treat first. The following proposition can be proved easily. It has been stated partly as Proposition 4.42 of \cite{M1}.  
\begin{prop}\label{AC1} For $J\subseteq I$ set ${\mathcal E}_J:=\Mklz{\ve{R}}{R\in \RkX_J}$. Then
\begin{eqnarray*}
  \WeD_J := \We_J \,{\mathcal E}_J =  {\mathcal E}_J \,\We_J
\end{eqnarray*}
is a submonoid of $\WeD$ canonically isomorphic to $\WeD(A_J)$. In particular, it holds $(\WeD_J)^\times=\WeD_J\cap\We=\We_J$. Furthermore, for $J,\,K\subseteq I$ the following properties are satisfied:
\begin{itemize}
\item[(a)] $\WeD_J\subseteq \WeD_K$ if and only if $J\subseteq K$.
\item[(b)] $\WeD_J\cap\WeD_K=\WeD_{J\cap K}$.
\item[(c)] $\WeD_J$ and $\WeD_K$ generate $\WeD_{J\cup K}$ as a monoid.
\end{itemize}
\end{prop}
We call $\WeD_J$ the {\it standard parabolic submonoid} of $\WeD$ of {\it type} $J$, and every conjugate $\sigma \WeD_J \sigma^{-1}$, $\sigma\in\We$, a {\it parabolic submonoid} of $\WeD$ of {\it type} $J$.\vspace*{1ex}

We take as Coxeter complex the set
\begin{eqnarray*}
   \mathcal C:=\Mklz{\sigma \We_J}{\sigma\in\We,\;J\subseteq I}
\end{eqnarray*}
partially ordered by the reverse inclusion, i.e., for $\sigma,\,\sigma'\in\We$ and $J,\,J'\subseteq I$, 
\begin{eqnarray*}
    \sigma \We_J\leq \sigma'\We_{J'} \;:\iff\; \sigma \We_J\supseteq \sigma'\We_{J'}. 
\end{eqnarray*}
The Weyl group $\We$ acts order preservingly on $\mathcal C$ by multiplication from the left. We now look for an extension to an action of the face monoid $\WeD$ on $\mathcal C$.
%
%%%%%%%%%%%%%%%%%%%%%%%%%%%%%%%%%%%%%%%%%%%%%%%%%%%%%%%%%%%%%%%%%%%%%%%%%%%%%%%%%%%%%%%%%%%%%%%%%%%%%%%%%%%%%%%%%%%%%%%%%%%%%%%%%%%
%Please add some additional space:
\vspace*{1ex}
%
%%%%%%%%%%%%%%%%%%%%%%%%%%%%%%%%%%%%%%%%%%%%%%%%%%%%%%%%%%%%%%%%%%%%%%%%%%%%%%%%%%%%%%%%%%%%%%%%%%%%%%%%%%%%%%%%%%%%%%%%%%%%%%%%%%%%
%

There is an obvious action of the face monoid $\WeD$ on the Tits cone $X$, extending the action of $\We$. It has been stated without proof, which we add here, as Proposition 2.27 of \cite{M1}. 
\begin{prop}\label{AC2}
The face monoid $\WeD$ acts faithfully on the Tits cone $X$ by 
\begin{eqnarray*}
  \left(\sigma\ve{R}\right) \la :=  \left\{\begin{array}{ccc}
                                       \sigma\la &\mb{ if }& \la\in R \\
                                         0 &\mb{ if }& \la\in X\setminus R 
                                      \end{array}\right.   &\quad,\quad& \sigma\ve{R}\in \WeD .
\end{eqnarray*} 
The parabolic submonoid $\WeD_J$ is the stabilizer of any element of the facet $F_J$, $J\subseteq I$.
\end{prop} 
 \Proof The maps on $X$ given by the elements of $\WeD$ are well defined: If $\sigma \ve{R},\,\ti{\sigma}\ve{\ti{R}}\in\WeD$ such that 
$\sigma \ve{R}=\ti{\sigma} \ve{\ti{R}}$ then $R=\ti{R}$ and $\sigma^{-1}\ti{\sigma}\in Z_{\mathcal W}(R)$. This implies immediately 
$\sigma \ve{R} \la = \ti{\sigma}\ve{\ti{R}}\la $ for all $\la\in X$.

Trivially, $1=\ve{X}$ acts as identity on $X$. Now let  $\sigma \ve{R},\,\tau \ve{S}\in\WeD$. For all $\la\in X$ we find by the definition of the action and by the definition of the multiplication of $\WeD$: 
\begin{eqnarray*}
   &&  \sigma\ve{R}(\tau \ve{S}\la) =  \left\{\begin{array}{cc}
                                       \sigma\tau\la &\mb{ if }\,\la\in S \,\mb{ and }\,\tau \la\in R\\
                                         0 &\mb{ else }
                                      \end{array}\right\} \\
   &&                               =  \left\{\begin{array}{ccc}
                                        \sigma\tau\la &\mb{ if }\,\la\in \tau^{-1}R\cap S\\
                                         0 &\mb{ else }
                                      \end{array}\right\} 
          =\sigma\tau \ve{\tau^{-1}R\cap S}\la =  (\,\sigma\ve{R}\tau \ve{S}\,)\,\la  .
\end{eqnarray*}

Now let $\sigma \ve{R},\,\tau \ve{S}\in\WeD$ such that $\sigma \ve{R}\la=\tau \ve{S}\la$ for all $\la\in X$. By the definition of the action this implies 
$R=S$ and $\sigma\la=\tau\la$ for all $\la\in R=S$, which by definition is equivalent to $\sigma \ve{R}=\tau \ve{S}$. 

At last we compute the stabilizer of the point $\la\in F_J$, $J\subseteq I$. Note that $F_J\subseteq ri(R(J^\infty)) \subseteq R(J^\infty)$ by Theorem \ref{FL2}. 
The stabilizer of $\la$ contains $\We_J$ because of $\la\in F_J$. The stabilizer of $\la$ contains ${\mathcal E}_J$ because of $\la\in R(J^\infty)$. 
Therefore, it also contains the parabolic submonoid $\WeD_J={\mathcal E}_J\We_J=\We_J {\mathcal E}_J$. 
Now let $\sigma \ve{R}\in\WeD$ such that $\sigma \ve{R}\la=\la$. By the definition of the action we find $\sigma\la=\la\in R$. It follows $\sigma\in\We_J$. 
Because of $\la\in ri(R(J^\infty))$, it also follows $R(J^\infty)\subseteq R$. Therefore, $\sigma\ve{R}\in\WeD_J$.
\qed 

It is very natural to consider this action because it is directly related to representation theory: Let $V$ be a module used for the Tannaka recponstruction of $\GD$ in \cite{M1}. In particular, $V$ is a direct sum of highest weight modules $L(\La)$, $\La\in P^+$. Let $\hat{n}_{\hat{\sigma}}\in \ND$ be an element which projects to $\hat{\sigma}\in \WeD=\ND/T$. It is easy to check that
\begin{eqnarray*}
    \hat{n}_{\hat{\sigma}} V_\la\subseteq  V_{\hat{\sigma} \la}\;\mb{ where }\; \la\in P(V)\cup \{0\},
\end{eqnarray*}
use Proposition \ref{FM4} and the definition of $e(R)$ in (\ref{defeR}), $R$ a face of the Tits cone.
%%%%%%%%%%%%%%%%%%%%%%%%%%%%%%%%%%%%%%%%%%%%%%%%%%%%%%%%%%%%%%%%%%%%%%%%%%%%%%%%%%%%%%%%%%%%%%%%%%%%%%%%%%%%%%%%%%%%%%%%%%%%%%%%%%%%
%Please add some additional space:
\vspace*{1ex}
%
%%%%%%%%%%%%%%%%%%%%%%%%%%%%%%%%%%%%%%%%%%%%%%%%%%%%%%%%%%%%%%%%%%%%%%%%%%%%%%%%%%%%%%%%%%%%%%%%%%%%%%%%%%%%%%%%%%%%%%%%%%%%%%%%%%%%
%

The set of open facets of the Tits cone $X$ identifies with the Coxeter complex $\mathcal C$, identifying the open facet $\tau F_J$ of the Tits cone with the element $\tau \We_J$ of the Coxeter complex for all $\tau\in\We$, $J\subseteq I$.
The partial order on open facets is the inclusion order of the corresponding closed facets. The action of the Weyl group $\We$ on 
the open facets is induced by the the action of $\We$ on the Tits cone $X$.

Now, the action of the face monoid $\WeD$ on the Tits cone $X$, which is given in Proposition \ref{AC2}, induces an action of $\WeD$ on the open facets. The general situation is described by the following Lemma, which is trivial to prove:
\begin{lem}\label{ACL} Let $M$ be a monoid acting on a nonempty set $S$. Let $\mathcal P$ be a partition of $S$ with the property: For every $m\in M$ and $P\in{\mathcal P}$ there 
exists a set $P(m)\in {\mathcal P}$ such that $m P\subseteq P(m)$. Then we get an action $\mb{}_\bullet$ of $M$ on $\mathcal P$ by
\begin{eqnarray*}
  m_\bullet P:=P(m)\quad \mb{ where }\quad  m\in M,\;P\in{\mathcal P}.
\end{eqnarray*}
\end{lem}

Specializing to the action of $\WeD$ on the Tits cone $X$ given in Proposition \ref{AC2} and to the partition of open facets, identifying the open facets with the elements of the Coxeter complex, we arrive at the following proposition. We use the notation 
\begin{eqnarray*}
   \We_J\downarrow :=\Mklz{V\in{\mathcal C}}{V\leq \We_J}, \quad J\subseteq I.
\end{eqnarray*}
\begin{prop}[Bad action]\label{AC3}
The face monoid $\WeD$ acts faithfully on the Coxeter complex $\,{\mathcal C}$ as a set by 
\begin{eqnarray}\label{baCox}
  \sigma_1\ve{R(\Th)}\sigma_2  \bad \tau \We_J :=  \left\{\begin{array}{ccc}
                                       \sigma_1\sigma_2\tau  \We_J&\mb{ if }&  \sigma_2\tau \We_J\in  \We_{\Th^\bot} (\We_\Th \downarrow )\\
                                         \We_I &\mb{ else }& 
                                      \end{array}\right. ,
\end{eqnarray} 
where $\sigma_1\ve{R(\Th)}\sigma_2\in\WeD$ and $\tau \We_J\in {\mathcal C}$. 
Unfortunately, this action does not preserve the order on $\mathcal C$. The parabolic submonoid $\WeD_J$ is the stabilizer of $\We_J\in {\mathcal C}$, $J\subseteq I$.
\end{prop} 
 \Proof Consider the action of $\WeD$ on the Tits cone $X$ given in Proposition \ref{AC2}. The faces of $X$ are unions of open facets of $X$. Therefore, by the 
definition of the action of $\WeD$, every element of $\WeD$ maps any open facet of $X$ into an open facet of $X$. By the last Lemma we get an action of $\WeD$ on 
the set of open facets. If we identify the open facets with the elements of the Coxeter complex, i.e., $\tau F_J$ identifies with $\tau\We_J$, 
where $\tau\in\We$ and $J\subseteq I$, then this action gives the action of the proposition. For simplicity, we work in the remaining proof with the action of $\WeD$ 
on the open facets.

Let $\la\in X$ and let $F$ be an open facet. We next show $|\WeD\la\cap F|\leq 1$: If $\la\in X\setminus F_I$ then $\WeD\la\subseteq \We\la\cup\{0\}$ and $\We\la\cap F_I=\emptyset$. It follows that $|\WeD\la\cap F_I|\leq 1$ and 
$|\WeD\la \cap F|\leq |\We\la\cap F|\leq 1$ for every facet $F\neq F_I$ of $X$. If $\la\in F_I$ then $\WeD\la =\We\la=\la$.

Now we compute the stabilizer of $F_J$, $J\subseteq I$. 
For the action of $\We$ given in Proposition \ref{AC2} the monoid $\WeD_J$ is the stabilizer of every point of $F_J$. Therefore, $\WeD_J$ is contained in the 
stabilizer of $F_J$ under the action $\bad$ of $\WeD$ on the set of open facets. 
Now let $\widehat{\sigma}\in\WeD$ such that $\widehat{\sigma}\bad F_J=F_J$. Choose an element $\la\in F_J$. Then $\widehat{\sigma}\la\in F_J$. Since 
$|\WeD\la\cap F_J|\leq 1$ we find $\widehat{\sigma} \la =\la$. By Proposition \ref{AC2} it follows $\widehat{\sigma}\in\WeD_J$.

Next we show that $\WeD$ acts faithfully on the Coxeter complex. Let  $\widehat{\sigma},\,\widehat{\tau}\in\WeD$ such that 
$\widehat{\sigma}\bad F=\widehat{\tau}\bad F$ for every open facet $F$ of $X$. 
Let $\la\in X$ and let $F$ be the open facet which contains $\la$. It holds $\widehat{\sigma}\la \in\widehat{\sigma}\bad F=\widehat{\tau}\bad F \ni \widehat{\tau}\la$. Since the orbit $\WeD\la$ intersects the facet $\widehat{\sigma}\bad F=\widehat{\tau}\bad F$ at most in one point we get $\widehat{\sigma}\la=\widehat{\tau}\la$. 
Since the action of $\WeD$ on the Tits cone $X$ given in Proposition \ref{AC2} is faithful, we get $\widehat{\sigma}=\widehat{\tau}$.

It remains to give an example that the action is in general not order preserving:   
Let $A$ be of indefinite, but not strongly hyperbolic type. Then there exists a nontrivial special set $\emptyset\neq \Th\neq I$. 
Now $F_\emptyset \gneqq F_\Th$, but $\ve{R(\Th)}\bad F_\emptyset = F_I \lneqq F_\Th=\ve{R(\Th)}\bad F_\Th$.   
\qed 

This Proposition is disappointing. The representation theoretically motivated obvious extension of the action of $\We$ on $\mathcal C$ to an action of $\WeD$ on $\mathcal C$ is not order preserving any more. As seen in the last proof, responsible for this breakdown is the trivial ``else''-part of formula (\ref{baCox}).

We want to keep the essential properties of the bad action, but change the ``else''-part of formula (\ref{baCox}) to obtain order preserving actions of $\WeD$ on $\mathcal C$. Which properties could we keep? These are the following properties:
\begin{itemize}
\item[(a)] The nontrivial "if"-part of formula (\ref{baCox}).
\item[(b)] The stabilizers $\mb{Stab}_{\widehat{\mathcal W}}(\We_J)=\WeD_J$, $J\subseteq I$. Then the parabolic submonoids of $\WeD$ play the same role for the action of $\WeD$ as the parabolic subgroups of $\We$ do for the action of $\We$.
\end{itemize}

The next proposition gives different equivalent descriptions of the ``if''-part of formula (\ref{baCox}):
\begin{prop}\label{ACneu} Let $\bullet$ be an action of the face monoid $\WeD$ on the Coxeter complex $\mathcal C$ which extends the action of the Weyl group $\We$ on $\mathcal C$. The following three statements are equivalent:
\begin{itemize}
\item[(i)] Let $\sigma_1\ve{R(\Th)}\sigma_2\in\WeD$ and  $\tau\We_J\in {\mathcal C}$. Then
\begin{eqnarray*}
 \sigma_1\ve{R(\Th)}\sigma_2 \mb{}_\bullet \tau \We_J = 
                                       \sigma_1\sigma_2\tau  \We_J &\mb{ if }&  \sigma_2\tau \We_J\in  \We_{\Th^\bot} (\We_\Th \downarrow ).
\end{eqnarray*}
\item[(ii)] $\ve{R(J^\infty)}_\bullet\We_J=\We_J$ for all $J\subseteq I$.
\item[(iii)] $\mb{Stab}_{\widehat{\mathcal W}}(\We_J)\supseteq\WeD_J$ for all $J\subseteq I$.
\end{itemize}
\end{prop}
 \Proof Trivially, $(ii)$ follows from $(i)$. Now suppose that (ii) holds. Fix $J\subseteq I$. For every special set $\Th\subseteq J$ we find
\begin{eqnarray*}
&& \ve{R(\Th)}_\bullet \We_J = \ve{R(\Th)}_\bullet\left(\ve{R(J^\infty)}_\bullet\We_J\right) = \left(\ve{R(\Th)}\ve{R(J^\infty)}\right)_\bullet\We_J\\
&&\quad\qquad = \ve{R(\Th)\cap R(J^\infty)}_\bullet\We_J = \ve{R(\Th\cup J^\infty)}_\bullet\We_J = \ve{R(J^\infty)}_\bullet\We_J= \We_J 
\end{eqnarray*}
Also, $\We_J\subseteq \mb{Stab}_{\widehat{\mathcal W}}(\We_J)$. Therefore, it follows (iii). 

Now suppose that (iii) holds. Let $\sigma_1\ve{R(\Th)}\sigma_2\in\WeD$ and  $\tau\We_J\in {\mathcal C}$ such that $\sigma_2\tau \We_J\in  \We_{\Th^\bot} (\We_\Th \downarrow )$. Then $J\supseteq \Th$ and $\sigma_2\tau\in \We_{\Th^\bot} \We_J$. Write $\sigma_2\tau$ in the form $\sigma_2\tau=ab$ with $a\in\We_{\Th^\bot}$ and $b\in\We_J$. Then
\begin{eqnarray*}
&& \sigma_1\ve{R(\Th)}\sigma_2 \mb{}_\bullet \tau \We_J = \sigma_1\ve{R(\Th)}_\bullet\sigma_2\tau\We_J  =\sigma_1\ve{R(\Th)}_\bullet a\We_J \\
&&  \quad\qquad=\sigma_1\ve{R(\Th)}a_\bullet \We_J  = \sigma_1 a\ve{R(\Th)}_\bullet \We_J = \sigma_1 a \We_J =\sigma_1\sigma_2\We_J.
\end{eqnarray*}
\qed 

By this proposition, the stabilizer condition (b) is sufficient for the "if"-part condition (a). 
It is reasonable to keep condition (b) for modifications of the bad action, because we will see in Theorem \ref{AB12}: If for the action of $\WeD$ on the Coxeter complex $\mathcal C$ there exists a compatible action of $\GD$ on the building $\Omega$, and if $\mb{Stab}_{\widehat{\mathcal W}}(\We_J)\supseteq\WeD_J$ holds for all $J\subseteq I$, then already $\mb{Stab}_{\widehat{\mathcal W}}(\We_J)=\WeD_J$ for all $J\subseteq I$. 
\begin{defn}\label{DefgoodW} An action of the face monoid $\WeD$ on the Coxeter complex $\mathcal C$, extending the action of the Weyl group $\We$ on $\mathcal C$, is called good if the following two conditions hold:
\begin{itemize}
\item[(1)]$\mb{Stab}_{\widehat{\mathcal W}}(\We_J)=\WeD_J$ for all $J\subseteq I$.
\item[(2)] $\WeD$ acts order preservingly on $\mathcal C$.
\end{itemize}
\end{defn}

We next give two good actions, which we call good action 1 and 2. In particular, good action 2 is quite complicated.
We explain how both actions can be found. They are related to the Tits cone and to a version of the Tits cone given by E. J. Looijenga. 

The set of closed facets of the Tits cone $X$ also identifies with the Coxeter complex $\mathcal C$, identifying the closed facet $\tau \overline{F}_J$ of the Tits cone with the element $\tau \We_J$ of the Coxeter complex for all $\tau\in\We$, $J\subseteq I$. The partial order is the inclusion order of the closed facets. The action of the Weyl group $\We$ on the closed facets is induced by the action of $\We$ on the Tits cone $X$.

Recall from Section \ref{section2} that the intersection of a face of $X$ and a closed facet of $X$ is again a closed facet of $X$. Now, the idea for the first modified action of the face monoid $\WeD$ on the Coxeter complex is to try the following easy formula
\begin{eqnarray}\label{defgood1facets}
   \sigma \ve{R}\goodone \overline{F}:= \sigma (R\cap \overline{F}) \quad \mb{ where } \quad \sigma\ve{R}\in\WeD, \;\overline{F} \mb{ a closed facet}.
\end{eqnarray} 
It works. To be able to compare it with the bad action and later with good action 2 we use Theorem \ref{FL6} to rewrite this formula: 
\begin{prop}[Good action 1]\label{ACgood1} We get a faithful, good action $\goodone$ of the face monoid $\WeD$ on the Coxeter complex $\mathcal{C}$ in the following way: Let $\sigma_1\ve{R(\Th)}\sigma_2\in\WeD$, and $\tau \We_J\in {\mathcal C}$. Decompose $\sigma_2\tau$ in the form
\begin{eqnarray*}
  \sigma_2\tau = a_1 a_2 y c
\end{eqnarray*}
where $a_1\in \We_{\Th^\bot}$, $a_2\in\We_{\Th}$, $y\in\mb{}^{\Th\cup \Th^\bot}\We^J$ and $c\in \We_J$. Define
\begin{eqnarray}\label{defgood1}
  \sigma_1\ve{R(\Th)}\sigma_2 \goodone  \tau\We_J := \sigma_1 a_1\We_{\Th\cup J\cup \,\red{y}}.
\end{eqnarray}
\end{prop}
 \Proof We first show that the map $\sigma\ve{R}\goodone$ given by formula (\ref{defgood1facets}) is well defined. Let $\sigma,\,\tau\in\We$, and $R,\,S$ faces of $X$ such that $\sigma\ve{R}=\tau \ve{S}$. Then by definition $R=S$ and $\tau^{-1}\sigma \in Z_{\mathcal W}(R)$. It follows
\begin{eqnarray*}
    \sigma \ve{R}\goodone \overline{F} = \sigma (R\cap \overline{F})= \tau (\tau^{-1}\sigma) (R\cap \overline{F})
                                     = \tau (R\cap \overline{F})= \tau (S\cap \overline{F})=\tau \ve{S}\goodone \overline{F}
\end{eqnarray*} 
for every closed facet $\overline{F}$ of $X$.

Furthermore, $\goodone$ is an action: Trivially, $1\goodone \overline{F} = X\cap \overline{F} =\overline{F} $ for all closed facets $\overline{F}$. Let $\sigma_1\ve{R},\,\tau\ve{S}\in\WeD$. Then
\begin{eqnarray*}
  &&  \sigma\ve{R}\goodone (\tau\ve{S}\goodone\overline{F} ) = \sigma\ve{R}\goodone \tau (S\cap\overline{F} ) = \sigma(R\cap \tau (S\cap\overline{F})) \\
  &&\qquad\quad= \sigma\tau (\tau^{-1}R\cap S\cap \overline{F}) = \sigma\tau\ve{\tau^{-1}R\cap S}\goodone\overline{F} = (\sigma\ve{R}\tau\ve{S})\goodone \overline{F}
\end{eqnarray*} 
for all closed facets $\overline{F}$.
Obviously, the action of $\WeD$ on the closed facets is order preserving. 

Next we use Theorem \ref{FL6} to rewrite formula (\ref{defgood1facets}) in formula (\ref{defgood1}). Recall that $\We_{\Th\cup\Th^\bot}$ stabilizes the set 
$R(\Th)$ as a whole. It holds
\begin{eqnarray*}
  &&   \sigma_1\ve{R(\Th)}\sigma_2 \goodone  \tau \overline{F}_J = \sigma_1 (R(\Th)\cap a_1 a_2 y c \overline{F}_J) 
      = \sigma_1 a_1 a_2  ( R(\Th)\cap y  \overline{F}_J ) \\
  &&  \qquad\quad= \sigma_1 a_1 a_2   \overline{F}_{\Th\cup J\cup\,\red{y}} 
      = \sigma_1 a_1\overline{F}_{\Th\cup J\cup \,\red{y}}.
\end{eqnarray*}

Now we compute the stabilizer of the facet $\overline{F}_J$. Let $\sigma\ve{R}\in\We$ such that 
\begin{eqnarray}\label{stclf1}
   \overline{F}_J=\sigma\ve{R}\goodone \overline{F}_J= \sigma (R\cap \overline{F}_J).
\end{eqnarray} 
Comparing the open facets of type $J$ on the left and on the right we get $F_J=\sigma F_J$ and $F_{J}\subseteq R$. Note that $F_{J}=F_{J^\infty\cup J^0}\subseteq ri\ R(J^\infty)$. It follows
\begin{eqnarray}\label{stclf2}
    \sigma\in\We_J \quad \mb{ and }\quad  R(J^\infty)\subseteq R,
\end{eqnarray}
which in turn is equivalent to $\sigma\ve{R}\in\WeD_J$, compare Proposition \ref{FL5} and the definition of $\WeD_J$ in Proposition \ref{AC1}. Obviously, from (\ref{stclf2}) follows (\ref{stclf1}).

At last, we show that the action is faithful. Suppose that $\sigma\ve{R},\,\tau\ve{S}\in\WeD$ such that $\sigma \ve{R}\goodone \overline{F} = \tau \ve{S}\goodone \overline{F}$ for all closed facets $\overline{F}$. Then $\sigma (R\cap \overline{F}) = \tau (S\cap \overline{F})$ for all closed facets $\overline{F}$. The face $R$ is a union of closed facets. Therefore, from this equation follows $\sigma R = \tau (S\cap R)$, from which in turn follows $\sigma R \subseteq  \tau S $ and $\sigma R \subseteq \tau R $. From the second inclusion we get $\tau^{-1}\sigma \in Z_{\mathcal W}(R)$, compare Theorem \ref{FL4} (a). By the first inclusion it follows $R\subseteq S$. Similarly, changing the roles of $R$ and $S$ we get $S\subseteq R$. Therefore, $\sigma\ve{R}=\tau\ve{S}$. 
\qed 

The bad action of $\WeD$ is induced by an action of $\WeD$ on the Tits cone $X$. Similarly, the second modified action of $\WeD$ is induced by an action of 
$\WeD$ on a space $X^\vee_{Loo}$ resembling a Tits cone, which has been introduced by E. Looijenga in \cite{Loo}. 
The definition of $X^\vee_{Loo}$ in \cite{Loo} is an ad-hoc definition. We give an interpretation of $X^\vee_{Loo}$ which also helps to find the action of $\WeD$ on $X^\vee_{Loo}$.  

Recall from Theorem 5.2 and Proposition 1.4 of \cite{M1} that the Zariski closure of the torus $T$ is a monoid $\TD$, which shares many properties with a toric 
monoid. In particular, it is a union of tori, i.e., 
\begin{eqnarray*}
  \TD=\dot{\bigcup_{R\;a\;face\;of\;X}} T e(R),
\end{eqnarray*}
where $T e(R)$ with the restricted multiplication is a torus with unit $e(R)$.  
It is easy to check that the monoid $\ND$ is an inverse monoid, the inverse map $\mb{}^{inv}:\ND\to\ND$ given by
\begin{eqnarray*}
     (n e(R))^{inv}= e(R) n^{-1} \quad \mb{ where }n\in N,\; R \mb{ a face of }X.
\end{eqnarray*}
Now we get an action $\mb{}_\bullet$ of $\ND$ on $\TD$ by
\begin{eqnarray*} 
  \hat{n}_\bullet \hat{t}:= \hat{n}\hat{t}\hat{n}^{inv} \;\mb{ where }\;\hat{n}\in\ND,\; \hat{t}\in \TD,
\end{eqnarray*}
which factors to an action $\mb{}_\bullet$ of $\WeD$ on $\TD$. Note that these actions extend the action $\mb{}_\bullet$ of $N$ on $T$ by conjugation and the 
corresponding action $\mb{}_\bullet$ of $\We$ on $T$. Explicitely, the action of $\WeD$ on $\TD$ is given by
\begin{eqnarray*}
 \sigma\ve{R}_\bullet t e(\ti{R}) = (\sigma_\bullet t) e(\sigma(R\cap \ti{R})) \;\mb{ where }\;\sigma\ve{R}\in\WeD,\; t e(\ti{R})\in\TD.
\end{eqnarray*}

Let $R=\sigma R(\Th)$. The torus $Te(R)$ and its (abelian) Lie algebra can be described up to isomorphisms as the following quotients:
\begin{eqnarray*}
   T e(R) \cong T/T_R \;&\mb{ where }&\; T_R:=\sigma_\bullet \prod_{i\in\Th} t_{h_i}(\F^\times),\\
     Lie(T e(R)) \cong \h/\h_R \;&\mb{ where }&\; \h_R:=\sigma \,\mb{span}\Mklz{h_i}{i\in\Th}.
\end{eqnarray*}
We now tie together the Lie algebras of these tori in a structure, which resembles an abelian Lie cone. 
\begin{prop}\label{AC4} For a face $R$ of $X$ denote by $p_R:\h\to\h/\h_R$ the canonical projection. 
On the set
\begin{eqnarray*}
 \widehat{\h}:=\dot{\bigcup_{R\;a\;face\;of\;X}} \h/\h_R
\end{eqnarray*}
there exist operations $+:\widehat{\h}\times \widehat{\h}\to \widehat{\h}$ and $\cdot:\F\times\widehat{\h}\to \widehat{\h}$ defined by
\begin{eqnarray*}
  p_R(h)+p_{\ti{R}}(\ti{h}):=p_{R\cap\ti{R}}(h+\ti{h}) \;\mb{ and  } \; a \cdot p_R(h):= p_R(a h), 
\end{eqnarray*}
where $R$, $\ti{R}$ are faces, $h,\,\ti{h}\in\h$, and $a\in\F$. Furthermore, $(\widehat{\h},\,+,\,\,\cdot\,)$ has similar properties as a linear space, except that 
$(\widehat{\h},\,+\,)$ is only an abelian monoid.

There exists an action of $\WeD$ on $\widehat{\h}$ defined by
\begin{eqnarray}\label{alc}
   \sigma\ve{R}_\bullet p_T (h) :=  p_{\sigma(R\cap T )} (\sigma h)\;\mb{ where }\;\sigma\ve{R}\in\WeD,\;T\;\mb{a face}, \;h\in\h.
\end{eqnarray}
\end{prop}
 \Proof The only statements which are not immediate are that the operations $+$ and $\mb{}_\bullet$ are well defined. 
We give here a direct algebraic proof. It would be also possible to obtain the operations $+$ and $\mb{}_\bullet$ by differentiating the corresponding operations on $\TD$. 

The operation $+$ is well defined if and only if 
\begin{eqnarray}\label{eqnhf}
  \h_R+\h_{\ti{R}}\subseteq \h_{R\cap \ti{R}}
\end{eqnarray}
for all faces $R$, $\ti{R}$. For any face $S$ the orthogonal complement of $\h_S$ in $\h^*$ is obtained by $(\h_S)^\bot =\mb{span}(S\cap P)$, compare Theorem \ref{FL2}. Taking the orthogonal complement, the inclusion (\ref{eqnhf}) is equivalent to 
\begin{eqnarray*}
     \mb{span}(R\cap P)\cap \mb{span}(\ti{R}\cap P)=     (\h_R)^\bot\cap (\h_{\ti{R}})^\bot \supseteq (\h_{R\cap \ti{R}})^\bot =\mb{span}((R\cap\ti{R})\cap P),
\end{eqnarray*} 
which holds trivially. 

To show that the operation $\mb{}_\bullet$ is well defined let $\sigma \ve{R}=\sigma'\ve{R'}$ and $p_T(h)=p_T(h')$. Then $R'=R$, $\sigma^{-1}\sigma'\in Z_{\mathcal W}(R)$, 
and $h-h'\in\h_T$. It follows
\begin{eqnarray}\label{Liecone1}
   \sigma(R\cap T)=\sigma (\sigma^{-1}\sigma') (R\cap T)= \sigma' (R\cap T)= \sigma' (R'\cap T). 
\end{eqnarray}
Let $R=\tau R(\Th)$. Because of $Z_{\mathcal W}(R)=\tau\We_\Th\tau^{-1}$ we have 
\begin{eqnarray*}
  h'-\eta h' \in \tau \,\mb{span}\Mklz{h_i}{i\in\Th}=\h_R \;\mb{ for all }\;\eta\in Z_{\mathcal W}(R).
\end{eqnarray*}
Therefore, by (\ref{eqnhf}),
\begin{eqnarray}\label{Liecone2}
   h-\sigma^{-1}\sigma'h'=  h-h'  +  h'-\sigma^{-1}\sigma'h' \in \h_T+\h_R\subseteq \h_{R\cap T}.  
\end{eqnarray}
From (\ref{Liecone1}) and (\ref{Liecone2}) follows
\begin{eqnarray*}
       \sigma \ve{R}_\bullet p_T (h) =  \sigma h +\h_{\sigma(R\cap T)} = \sigma' h' +\h_{\sigma(R\cap T)}  
      = \sigma' h' +\h_{\sigma'(R'\cap T)} =\sigma' \ve{R'}_\bullet p_T (h').
\end{eqnarray*}
\qed 

If we work over the field $\F=\R$ then, as usual, we equip the spaces $\widehat{\h}$, $\h$ by the lower index $\R$. To make our notation not too 
complicated, we don't equip the maps $p_R$, $R$ a face, with the index $\R$. Recall that $X^\vee$ denotes the Tits cone in $\h_\R$, and $(X^\vee)^\circ$ 
its interior. Now, it is easy to see that
\begin{eqnarray*}
   X_{Loo}^\vee:=\dot{\bigcup_{R\; a\; face}} p_R((X^\vee)^\circ)\;\subseteq \;\widehat{\h}_\R
\end{eqnarray*}
is a $\WeD$-invariant convex cone in $\widehat{\h}_\R$. As a $\We$-space, $X_{Loo}^\vee$ has been introduced in an ad hoc way by E. Looijenga in \cite{Loo}.\vspace*{1ex}

E. Looijenga also defined a $\We$-invariant partition of $X_{Loo}^\vee$. It can be obtained as follows:
Let $R(\Th)$ be a special face of $X$. The linear dual $(\h_\R /(\h_\R)_{R(\Th)})^*$ identifies in the obvious way with
\begin{eqnarray*}
  \mb{span}(R(\Th))=\Mklz{\la\in\h_\R^*}{\la(h_i)=0 \mb{ for all }i\in\Th} .
\end{eqnarray*}  
Furthermore,the families
\begin{eqnarray*}
  (h_i+ (\h_\R)_{R(\Th)}  )_{i\in\Th^\bot }\subseteq \h_\R/(\h_\R)_{R(\Th)} \;\mb{ and }\;
  (\al_i)_{i\in\Th^\bot}\subseteq \mb{span}(R(\Th))
\end{eqnarray*}
are a free realization over $\R$ for the generalized Cartan submatrix $A_{\Th^\bot}$. By Lemma 2.8 of \cite{Loo} this induces a partition of 
$p_{R(\Th)} ((X^\vee)^\circ)$ into the facets
\begin{eqnarray*}
      \tau F_J^\vee (\Th) \;\mb{ where }\; \tau\in\We_{\Th^\bot},\;J\subseteq\Th^\bot \mb{ such that }J=J^0,
\end{eqnarray*} 
and 
\begin{eqnarray*}
      F_J^\vee (\Th):= \left\{h+(\h_\R)_{R(\Th)} \in p_{R(\Th)} ((X^\vee)^\circ)\left|
           \,\begin{array}{l}\al_i(h)=0 \mb{ for }i\in J,\\ \al_i(h)>0 \mb{ for }i\in \Th^\bot\setminus J.\end{array}
\right.\right\}.
\end{eqnarray*} 
Furthermore, also by Lemma 2.2 of \cite{Loo}, the map $p_{R(\Th)}$ maps every facet of $(X_{Loo}^\vee)^\circ$ into a facet of $p_{R(\Th)}((X_{Loo}^\vee)^\circ)$. In addition, every facet of $p_{R(\Th)}((X_{Loo}^\vee)^\circ)$ is obtained as image of a unique facet of $Star(F_{\Th^\bot}^\vee)\cap (X^\vee)^\circ$. In the following theorem we describe this map on the facets explicitely. 
\begin{thm}\label{AC5} Let $\Th$ be special. Let $\tau F_J^\vee\subseteq (X^\vee)^\circ$, i.e., $\tau\in\We$ and $J=J^0$.  
Decompose $\tau$ in the form $\tau= a y c$ where $a\in \We_{\Th^\bot}$, $y\in\mb{}^{\Th^\bot}\We^J$, and $c\in \We_J$. Then
\begin{eqnarray}\label{Fvee}
    p_{R(\Th)}(\tau F_J^\vee)\subseteq a F_{\Th^\bot \cap y J}^\vee (\Th).
\end{eqnarray}
Equality holds for $\tau F_J^\vee\in Star(F_{\Th^\bot}^\vee)\cap (X^\vee)^\circ$. 
\end{thm}
 \Proof Let $a\in \We_{\Th^\bot}$ and $h\in (X^\vee)^\circ$. Then  
\begin{eqnarray}\label{Fvee1}
   p_{R(\Th)} (a h)= a\,  p_{R(\Th)} (h)
\end{eqnarray}
holds trivially.

Now let $h\in F_J^\vee$ and $y\in \mb{}^{\Th^\bot}\We^J$. Let $i\in\Th^\bot$. Because of $y^{-1}\al_i \in \pW$ we get
\begin{eqnarray*}
  \al_i( y h)= (y^{-1}\al_i) (h)\geq 0 \;\mb{ for all }\; i\in\Th^\bot.
\end{eqnarray*}
In addition, $\al_i( y h)= (y^{-1}\al_i) (h)=0$ if and only if $y^{-1}\al_i\in \We_J J$. Therefore, we have shown
\begin{eqnarray}\label{Fvee2}
    p_{R(\Th)} (y F_J^\vee)\subseteq F^\vee_{y{\mathcal W}_J J\cap \Th^\bot}(\Th).
\end{eqnarray}
By Theorem 2.74 of \cite{Ca}, which holds also for general Coxeter groups, $y{\mathcal W}_J y^{-1}\cap \We_{\Th^\bot}=\We_{y J\cap \Th^\bot}$. Therefore, $y\We_J J\cap \Th^\bot = y J\cap \Th^\bot$.

The element $c\in\We_J$ stabilizes $F_J^\vee$. From equation (\ref{Fvee1}) and inclusion (\ref{Fvee2}) follows inclusion (\ref{Fvee}) immediately. 
The statement about the equality is valid by Lemma 2.2 of \cite{Loo}.
\qed

Now by the $\We$-action on $X_{Loo}^\vee$, the partitions of $p_{R(\Th)} ((X^\vee)^\circ)$, $\Th$ special, extend to a $\We$-invariant partition of $X_{Loo}^\vee$ into facets.

We summarize the properties of the $\WeD$-action on $X_{Loo}^\vee$:
\begin{prop}\label{AC6}$\WeD$ acts faithfully on $X_{Loo}^\vee$ by
\begin{eqnarray*}
\sigma \ve{R} p_{\ti{R}}(h)=p_{\sigma(R\cap\ti{R})}(\sigma h) \;\mb{ where }\; h\in(X^\vee)^\circ.
\end{eqnarray*}
The parabolic submonoid $\WeD_J$ is the stabilizer of any element of the facet $F_{J^0}^\vee(J^\infty)$, $J\subseteq I$.
\end{prop}
 \Proof It remains to show that $\WeD$ acts faithfully and to compute the stabilizers.

Let $\sigma \ve{R} p_T(h)=\tau\ve{S}p_T(h)$ for all $h\in(X^\vee)^\circ$, and faces $T$. In particular, for $T=X$ we get
$p_{\sigma R } (\sigma h)=p_{\tau S}(\tau h)$ for all $h\in(X^\vee)^\circ$. This is equivalent to $\sigma R=\tau S$ and 
$\sigma h -\tau h\in(\h_\R)_{\sigma R}$ for all $h\in(X^\vee)^\circ$, which is in turn equivalent to
\begin{eqnarray*}
   S=\tau^{-1}\sigma R \;\mb{ and }\; h -\sigma^{-1}\tau h\in (\h_\R)_R \;\mb{ for all }\; h\in(X^\vee)^\circ.
\end{eqnarray*}
Because of $\left((\h_\R)_R\right)^\bot = \mb{span}(R)$ we find
\begin{eqnarray*}
  \tau^{-1}\sigma \phi(h) =\phi(\sigma^{-1}\tau h)=\phi(h) \;\mb{ for all }\; \phi\in \mb{span}(R),\; h\in(X^\vee)^\circ.
\end{eqnarray*}
Since the cone $(X^\vee)^\circ$ spans $\h_\R$, this equation holds for all $h\in \h_\R$. It follows $\tau^{-1}\sigma\in Z_{\mathcal W}(R)$. 
We get $S=\tau^{-1}\sigma R=R$. Both implies $\sigma \ve{R} =\tau\ve{S}$.

Let $h\in (X^\vee)^\circ$ such that $p_{R(J^\infty)}(h)\in F_{J^0}^\vee(J^\infty)$. For an element $\sigma \ve{R}\in\WeD$ the stabilizer equation
\begin{eqnarray*}
  p_{R(J^\infty)}(h)=\sigma \ve{R} p_{R(J^\infty)}(h)= p_{\sigma(R\cap R(J^\infty))}(\sigma h)
\end{eqnarray*}
is equivalent to
\begin{eqnarray}\label{stabeqn}
\sigma(R\cap R(J^\infty))=R(J^\infty) \;\mb{ and }\; \sigma h-h\in (\h_\R)_{R(J^\infty)}.
\end{eqnarray}
By the first equation it follows $ R(J^\infty)\supseteq \sigma^{-1}R(J^\infty)$, which is equivalent to $R(J^\infty)=\sigma^{-1}R(J^\infty)$, which in turn is equivalent to 
$\sigma\in\We_{J^\infty\cup (J^\infty)^\bot}$. Furthermore, it follows $R\supseteq R(J^\infty)$.
Write $\sigma$ in the form $\sigma=\sigma_\bot\sigma_0$ with $\sigma_\bot\in\We_{(J^\infty)^\bot}$ and $\sigma_0\in \We_{J^\infty}$. Because of 
$\sigma_0 h-h\in (\h_\R)_{R(J^\infty)}$ and $\sigma_\bot (\h_\R)_{R(J^\infty)}=(\h_\R)_{R(J^\infty)}$, the condition
\begin{eqnarray*}
  \sigma_\bot (\sigma_0 h -h) + \sigma_\bot h - h =\sigma h-h\in (\h_{R(J^\infty)})_\R
\end{eqnarray*}
is equivalent to $\sigma_\bot h-h\in (\h_\R)_{R(J^\infty)}$, which is equivalent to $\sigma_\bot p_{R(J^\infty)}(h)=p_{R(J^\infty)}(h)$. Since $p_{R(J^\infty)}(h)\in F_{J^0}^\vee(J^\infty)$, we get $\sigma_\bot \in\We_{J^0}$ and $\sigma=\sigma_\bot\sigma_0 \in\We_{J^0}\We_{J^\infty}=\We_J$. It follows $\sigma \ve{R}\in\WeD_J$. 
The other way round, it is easy to check that the conditions (\ref{stabeqn}) are satisfied for every element $\sigma\ve{R}\in\WeD_J$, i.e., 
$\sigma\in\We_J$ and $R\supseteq R(J^\infty)$. 
\qed 

To proceed further, we need the following Proposition:
\begin{prop}Let $\sigma\in\We$. Let $J,\, K \subseteq I$. Then $\sigma$ can be uniquely written in the form
\begin{eqnarray*}
  \sigma = a x c_1 c_2
\end{eqnarray*}
with $a\in\We_K$, $x\in\mb{}^K\We^{J^\infty\cup (J^\infty)^\bot}$, $c_1\in(\We_{(J^\infty)^\bot})^{J^0}$, $c_2\in\We_J$ such that $xc_1\in\mb{}^K\We^{J}$ and 
$xc_1c_2\in\mb{}^K\We$.
\end{prop}
 \Proof The uniqueness of the decomposition follows easily, because for $L\subseteq I$ and $\tau\in\We$ the decompositions $\tau=\tau_1\tau_2$ with $\tau_1\in\We_L$, $\tau_2\in\mb{}^L\We$ resp. $\tau_1\in\We^L$, $\tau_2\in\We_L$ are unique.

By Proposition 2.7.5. of \cite{Ca}, which holds also for infinite Coxeter groups, we may decompose $\sigma$ in the form $\sigma=axc$ with $a\in\We_K$, $x\in\mb{}^K\We^{J^\infty\cup(J^\infty)^\bot}$ and $c\in\We_{J^\infty\cup(J^\infty)^\bot}\cap \mb{}^{(J^\infty\cup(J^\infty)^\bot)\cap x^{-1}K}\We$.

To prove $xc\in\mb{}^K\We$ we show $(xc)^{-1}\al_i=c^{-1}x^{-1}\al_i\in\pW$ for all $i\in K$. Because of $x\in\mb{}^K\We$ we have $x^{-1}\al_i\in \pW\cap x^{-1}K$ for all $i\in K$. First let $i\in K$ such that $x^{-1}\al_i\in \pW\setminus \W^+_{J^\infty\cup(J^\infty)^\bot}$. Because of $c^{-1}\in\We_{J^\infty\cup(J^\infty)^\bot}$ it follows $c^{-1}x^{-1}\al_i\in\pW$. Now let $i\in K$ such that $x^{-1}\al_i\in \W^+_{J^\infty\cup(J^\infty)^\bot}$. Write $x^{-1}\al_i$ in the form
\begin{eqnarray*}
 x^{-1}\al_i=\sum_{j\in J^\infty\cup (J^\infty)^\bot}n_j\al_j 
\end{eqnarray*}
where $n_j\in\Nn$ for all $j\in J^\infty\cup (J^\infty)^\bot$. Then
\begin{eqnarray*}
     \al_i=\sum_{j\in J^\infty\cup (J^\infty)^\bot}n_j x\al_j. 
\end{eqnarray*}
Because of $x\in\We^{J^\infty\cup (J^\infty)^\bot}$ we have $x\al_j\in\pW$ for all $j\in J^\infty\cup (J^\infty)^\bot$. Therefore, there exists an element $j\in J^\infty\cup (J^\infty)^\bot$ such that $\al_i=x\al_j$. In particular, $x^{-1}\al_i \in (J^\infty\cup (J^\infty)^\bot)\cap x^{-1}K$. Because of $c\in \mb{}^{(J^\infty\cup(J^\infty)^\bot)\cap x^{-1}K}\We$ it follows $c^{-1}x^{-1}\al_i\in\pW$.

Now decompose $c$ in the form $c=c_1\ti{c}_2\ti{c}_3$ with $c_1\in(\We_{(J^\infty)^\bot})^{J^0}$, $\ti{c}_2\in\We_{J^0}$, and $\ti{c}_3\in \We_{J^\infty}$. Trivially, $c_2:=\ti{c}_2\ti{c}_3 \in\We_J$. It remains to show $x c_1\in\mb{}^K\We^J$. We first show $x c_1\in\We^J$. Note that $x\in \We^{J^\infty\cup(J^\infty)^\bot}$ and $c_1\in (\We_{(J^\infty)^\bot})^{J^0}$. If $j\in J^\infty$, then $xc_1\al_j=x\al_j\in\pW$. If $j\in J^0$, then $xc_1\al_j\in x\W^+_{(J^\infty)^\bot}\subseteq\pW$. Now we show $xc_1\in \mb{}^K \We$. Let $i\in K$. Since $xc=x c_1c_2\in\mb{}^K\We$ we have 
\begin{eqnarray}\label{xc1c2}
(xc_1c_2)^{-1}\al_i=(c_2)^{-1} (c_1)^{-1} x^{-1}\al_i\in\pW.
\end{eqnarray}
Here, because of $x\in\mb{}^K\We$, $x^{-1}\al_i\in\pW\cap x^{-1}K$. Now we consider two cases. First let $(c_1)^{-1} x^{-1}\al_i\notin\W_{(J^\infty)^\bot}$. Since $(c_2)^{-1}\in\We_{(J^\infty)^\bot}$, from (\ref{xc1c2}) follows $(c_1)^{-1} x^{-1}\al_i\in \W^+\setminus\W^+_{(J^\infty)^\bot}$. Now let $(c_1)^{-1} x^{-1}\al_i\in\W_{(J^\infty)^\bot}$. Suppose that $(c_1)^{-1} x^{-1}\al_i\in\W^-_{(J^\infty)^\bot}$. Then 
\begin{eqnarray*}
  (c_1)^{-1}x^{-1}\al_i=\sum_{j \in (J^\infty)^\bot\setminus J^0} n_j\al_j +\sum_{j\in J^0} n_j \al_j
\end{eqnarray*}
with $n_j\in-\Nn$ for all $j\in (J^\infty)^\bot$. Since $c_2\in \We_J=\We_{J^\infty\cup J^0}$, from (\ref{xc1c2}) follows
\begin{eqnarray*}
  (c_1)^{-1}x^{-1}\al_i= \sum_{j\in J^0} n_j \al_j .
\end{eqnarray*}
Multiply this equation on the left and on the right with $c_1$. Since $x\in \mb{}^K\We $ and $c_1\in (\We_{(J^\infty)^\bot})^{J^0}$ it follows
\begin{eqnarray*}
       \underbrace{x^{-1}\al_i}_{\in \Delta^+}= \sum_{j\in J^0} \underbrace{n_j}_{\leq 0} \underbrace{c_1 \al_j}_{\in\Delta^+} .
\end{eqnarray*}
This implies $n_j=0$ for all $j\in J^0$, which contradicts $x^{-1}\al_i\neq 0$.
\qed 

The set of facets of $X_{Loo}^\vee$ identifies with the Coxeter complex $\mathcal C$ as a $\We$-set, identifying $\tau F_{J^0}(J^\infty)$ with $\tau \We_J$ for all $\tau\in\We$, $J\subseteq I$. Now the action of $\WeD$ on $X_{Loo}^\vee$ induces an action on the Coxeter complex, which is described by the following theorem.
\begin{thm}[Good action 2]\label{AC7} We get a faithful, good action $\goodtwo$ of the face monoid $\WeD$ on the Coxeter complex $\mathcal{C}$ in the following way: Let $\sigma_1\ve{R(\Th)}\sigma_2\in\WeD$, and $\tau \We_J\in {\mathcal C}$. Decompose $\sigma_2\tau$ in the form
\begin{eqnarray*}
  \sigma_2\tau = a_1 a_2 x c_1 c_2
\end{eqnarray*}
where $a_1\in \We_{\Th^\bot}$, $a_2\in\We_{\Th}$, $x\in\mb{}^{\Th\cup \Th^\bot}\We^{J^\infty\cup (J^\infty)^\bot}$, $c_1\in (\We_{(J^\infty)^\bot})^{J^0}$ and  
$c_2\in \We_J$, such that in addition $x c_1\in \mb{}^{\Th\cup \Th^\bot}\We^J$. Set $\Xi:=\Th\cup J^\infty \cup \red{x}$. Then 
\begin{eqnarray*}
    \sigma_1\ve{R(\Th)}\sigma_2 \goodtwo \tau \We_J = \sigma_1 a_1\We_{\Xi\cup (\Xi^\bot\cap c_1 J^0)}.
\end{eqnarray*}  
\end{thm}
\begin{rem} In general, good action 1 of $\WeD$ is different from good action 2 of $\WeD$, which can be seen as follows. Let $A$ be of 
indefinite but not strongly hyperbolic type. Then there exists a special set $\Th$ such that $\emptyset\neq\Th\neq I$. Since $A$ is indecomposable 
it holds $\Th\cup\Th^\bot\neq I$. Let $i\in I\setminus(\Th\cup\Th^\bot)$ and denote by $\sigma_i$ the corresponding simple reflection. Then 
$\ve{R(\Th)}\goodone\sigma_i\We_\emptyset=\We_{\Th\cup \{i\}}$ but $\ve{R(\Th)}\goodtwo\sigma_i\We_\emptyset=\We_\Th$.
\end{rem}
 \Proof  Consider the action of $\WeD$ on $X_{Loo}^\vee$. We first show that every element of $\WeD$ maps any facet of $X_{Loo}^\vee$ into a facet of $X_{Loo}^\vee$. Let $\sigma_1\ve{R(\Th)}\sigma_2 \in\WeD$ and $\tau F_{J^0}^\vee(J^\infty)$ be a facet of $X_{Loo}^\vee$. Decompose $\sigma_2\tau$ as in the theorem. Now $\Th\subseteq \Xi$ implies $\Th^\bot \supseteq \Xi^\bot$, from which follows $xc_1\in \mb{}^{\Xi^\bot}\We^{J^0}$. By the definition of the action of $\WeD$ on $X_{Loo}^\vee$, by Theorem \ref{AC5}, Proposition \ref{AC6}, and Theorem \ref{FL4} we find
\begin{eqnarray*}
     &&   \left(\sigma_1\ve{R(\Th)}\sigma_2 \right)\tau F_{J^0}^\vee(J^\infty)  = \sigma_1\ve{R(\Th)}a_1 a_2 x c_1 c_2 p_{R(J^\infty)}(F_{J^0}^\vee)\\
     &&  = \sigma_1 a_1 \ve{R(\Th)} x c_1 p_{R(J^\infty)}(F_{J^0}^\vee)  = \sigma_1 a_1 \ve{R(\Th)} p_{x c_1  R(J^\infty)}( x c_1  F_{J^0}^\vee)\\
     && = \sigma_1 a_1  p_{R(\Th)\cap x R(J^\infty)}( x c_1 F_{J^0}^\vee)  = \sigma_1 a_1  p_{R(\Xi)}( x c_1 F_{J^0}^\vee) 
      \subseteq   \sigma_1 a_1 F_{\Xi^\bot\cap x c_1 J^0}^\vee(\Xi).  
\end{eqnarray*}
Because of $\Xi^\bot\subseteq \red{x}^\bot$ it follows
\begin{eqnarray}\label{alternativ}
  \Xi^\bot\cap x c_1 J^0=x(x^{-1}\Xi^\bot\cap c_1 J^0)=x(\Xi^\bot\cap c_1 J^0)= \Xi^\bot\cap c_1 J^0.
\end{eqnarray}
By Lemma \ref{ACL} the monoid $\WeD$ acts on the set of facets of  $X_{Loo}^\vee$. Identify the facets of $X_{Loo}^\vee$ with the elements of the Coxeter complex $\mathcal C$, i.e., identify $\tau F_{J^0}(J^\infty)$ with $\tau \We_J$ for all $\tau\in\We$, $J\subseteq I$. Then the $\WeD$-action on the set of facets of $X_{Loo}^\vee$ coincides with $\WeD$-action on $\mathcal C$ of the theorem.

Next we show that $\WeD_J$ is the stabilizer $\We_J$. By our identification, Proposition \ref{AC6} shows $\WeD_J \goodtwo \We_J=\We_J$. Now let $\hat{\sigma}\in\WeD$ such that $\hat{\sigma}\goodtwo \We_J = \We_J$. Write $\hat{\sigma}$ in normal form I, i.e., $\hat{\sigma}=\sigma_1\ve{R(\Th)}\sigma_2$ with $\sigma_1\in\We^{\Th}$ and $\sigma_2\in\mb{}^{\Th\cup\Th^\bot}\We$, $\Th$ special. Decompose $\sigma_2=xc_1c_2$ with $x\in\mb{}^{\Th\cup\Th^\bot}\We^{J^\infty\cup (J^\infty)^\bot}$, $c_1\in (\We_{(J^\infty)^\bot})^{J^0}$, and $c_2\in\We_J$, such that $xc_1\in \mb{}^{\Th\cup\Th^\bot}\We^J$. Set $\Xi:= \Th\cup J^\infty\cup\red{x}$. Then 
\begin{eqnarray*}
   \We_J=\hat{\sigma}_\goodtwo\We_J =\sigma_1 \We_{\Xi\cup(\Xi^\bot\cap c_1 J^0)}.
\end{eqnarray*}
As we have seen above, $\Xi=\Th\cup J^\infty\cup\red{x}$ is a special set. Furthermore, $\Xi^\bot \cap c_1 J^0 $ has only components of finite type. It follows
\begin{eqnarray*}
    J^\infty=\Xi=\Th\cup J^\infty\cup\red{x}   \quad \mb{ and }\quad  J^0=\Xi^\bot \cap c_1 J^0 \quad\mb{ and }\quad \sigma_1 \in\We_J.
\end{eqnarray*}
Now, $\red{x}\subseteq J^\infty$ and $\Th\subseteq J^\infty$ imply $x\in\We_{J^\infty}\cap \mb{}^{\Th\cup\Th^\bot}\We^{J^\infty\cup (J^\infty)^\bot}=\{1\}$ and $\Th^\bot\supseteq (J^\infty)^\bot$. It follows $\We_{\Th\cup \Th^\bot}\sigma_2=\We_{\Th\cup \Th^\bot} c_2$. Since $\sigma_2$ is a minimal coset representative it holds $l(\sigma_2)\leq l(c_2)$. On the other hand also $l(\sigma_2)=l(x)+l(c_1)+l(c_2)=0+l(c_1)+l(c_2)$. It follows $c_1=1$ and $\sigma_2=c_2\in\We_J$. Therefore $\hat{\sigma}=\sigma_1\ve{R(\Th)}\sigma_2\in\WeD_J$.  

Next we show that the action of $\WeD$ is faithful. Let $\sigma_1\ve{R(\Th)}\sigma_2\in\WeD$ and $\sigma'_1\ve{R(\Th')}\sigma'_2\in\WeD$ such that 
\begin{eqnarray}\label{fgl}
  (\sigma_1\ve{R(\Th)}\sigma_2 )_\goodtwo\tau \We_J=(\sigma'_1\ve{R(\Th')}\sigma'_2 )_\goodtwo\tau \We_J  \quad \mb{ for all }\quad \tau \We_J\in {\mathcal C}.
\end{eqnarray}
In particular, for $J=\emptyset$ and $\tau=(\sigma'_2)^{-1}$ we get
\begin{eqnarray}\label{eqac1}
    \sigma_1\ve{R(\Th)}\sigma_2 (\sigma'_2)^{-1}_\goodtwo \We_\emptyset= \sigma'_1\ve{R(\Th')}_\goodtwo\We_\emptyset. 
\end{eqnarray}
Decompose $\sigma_2 (\sigma'_2)^{-1}$ in the form $\sigma_2 (\sigma'_2)^{-1}=a_1 a_2 x c_1 c_2$, with $a_1\in\We_{\Th^\bot}$, $a_2\in\We_{\Th}$, $x=1$, 
$c_1\in \mb{}^{\Th\cup \Th^\bot}\We$, and $c_2=1$. Then equation (\ref{eqac1}) is equivalent to $\sigma_1 a_1\We_{\Th}=\sigma'_1\We_{\Th'}$. Therefore,
\begin{eqnarray}\label{tsa}
   \Th'=\Th \quad\mb{ and }\quad (\sigma'_1)^{-1} \sigma_1 a_1\in \We_{\Th}.
\end{eqnarray}
Inserting $J=\Theta$ and $\tau=(\sigma'_2)^{-1}$ in equation (\ref{fgl}) we get
\begin{eqnarray}\label{eqac2}
     \sigma_1\ve{R(\Th)}\sigma_2 (\sigma'_2)^{-1}_\goodtwo \We_\Th= \sigma'_1\ve{R(\Th')}_\goodtwo\We_\Th.
\end{eqnarray}
Decompose $\sigma_2 (\sigma'_2)^{-1}$ in the form $\sigma_2 (\sigma'_2)^{-1}=a'_1 a'_2 x' c'_1 c'_2$ with $a'_1\in\We_{\Th^\bot}$, $a'_2\in\We_\Th$, $x'\in\mb{}^{\Th\cup \Th^\bot}\We^{\Th\cup \Th^\bot}$, $c'_1\in\We_{\Th^\bot}$, $c'_2\in\We_{\Th}$, such that 
$x'c'_1\in\mb{}^{\Th\cup \Th^\bot}\We^\Th$. 
Inserting $\Th'=\Th$, (\ref{eqac2}) is equivalent to $\sigma_1 a'_1\We_{\Th\cup\red{x'}}=\sigma'_1\We_\Th$. 
It follows $\red{x'}\subseteq \Th$, resp. $x'\in \We_{\Th}\cap \mb{}^{\Th\cup \Th^\bot}\We^{\Th\cup \Th^\bot}=\{1\}$.  
Therefore $c'_1= x'c'_1\in \mb{}^{\Th\cup \Th^\bot}\We^\Th\cap \We_{\Th^\bot}=\{1\}$. Since
$\sigma_2 (\sigma'_2)^{-1}=a_1a_2 c_1= a'_1 a'_2 c'_2$
%\begin{eqnarray*}
%  \sigma_2 (\sigma'_2)^{-1}=a_1a_2 c_1= a'_1 a'_2 c'_2
%\end{eqnarray*}
with $a_1\in\We_{\Th^\bot}$, $a_2\in\We_{\Th}$, $c_1\in \mb{}^{\Th\cup \Th^\bot}\We$, and $a'_1\in\We_{\Th^\bot}$, $a'_2\in\We_\Th$,  $c'_2\in\We_{\Th}$ 
we find  $a_1=a'_1$ and $a_2=a'_2 c'_2$ and $c_1=1$. With (\ref{tsa}) it follows
\begin{eqnarray*}
   && \sigma_1\ve{R(\Th)}\sigma_2 = \sigma_1\ve{R(\Th)}\sigma_2 (\sigma'_2)^{-1} \sigma'_2= \sigma_1 a_1 \ve{R(\Th)}\sigma'_2 \\
   && =\sigma'_1(\sigma'_1)^{-1} \sigma_1 a_1 \ve{R(\Th)}\sigma'_2 =\sigma'_1\ve{R(\Th)}\sigma'_2 = \sigma'_1\ve{R(\Th')}\sigma'_2 .
\end{eqnarray*}

Since the action of $\We$ is order preserving, it remains to show that any element $\ve{R(\Th)}$, $\Th$ special, acts order preservingly. Recall formula (\ref{alternativ}). Let $z_J \We_J,\,z_K\We_K\in{\mathcal C}$. The relation $z_J\We_J\leq z_K\We_K$ is equivalent to $K\subseteq J$ and $z_J^{-1}z_K\in\We_J$. Because of $z_J\We_J=z_K(z_K^{-1}z_J)\We_J=z_K\We_J$ it is sufficient to consider elements $z \We_J,\,z\We_K\in{\mathcal C}$ with $z\We_J\leq z\We_K$, which is equivalent to $K\subseteq J$. Trivially, it follows
\begin{eqnarray}\label{KJc}
   K^\infty\subseteq J^\infty \quad \mb{ and }\quad K^0 = \underbrace{(K^0\cap J^\infty)}_{\subseteq J^\infty} \cup \underbrace{(K^0\cap J^0)}_{\subseteq (J^\infty)^\bot}.
\end{eqnarray}

Decompose $z$ in the form $z= a_1 a_2 x c_1 c_2$ where $a_1\in \We_{\Th^\bot}$, $a_2\in\We_\Th$, 
$x\in\mb{}^{\Th\cup\Th^\bot}\We^{J^\infty \cup (J^\infty)^\bot}$, $c_1\in (\We_{(J^\infty)^\bot})^{J^0}$ 
and $c_2\in \We_J$, such that in addition $x c_1\in \mb{}^{\Th\cup\Th^\bot}\We^J$ and $xc_1c_2\in \mb{}^{\Th \cup \Th^\bot}\We$. Set 
$\Xi:=\Th\cup J^\infty \cup \red{x}$. Then
\begin{eqnarray*}
   \ve{R(\Th)} z\We_J = a_1\We_{\Xi\cup (\Xi^\bot\cap x c_1 J^0)}.
\end{eqnarray*}
Decompose $z$ in the form $z= a'_1 a'_2  x' c'_1 c'_2$ where $a'_1\in \We_{\Th^\bot}$, $a'_2\in\We_{\Th}$, 
$x'\in\mb{}^{\Th\cup \Th^\bot}\We^{K^\infty\cup (K^\infty)^\bot}$, 
$c'_1\in (\We_{(K^\infty)^\bot})^{K^0}$ and $c'_2\in \We_K$, such that in addition $ x' c'_1\in \mb{}^{\Th\cup \Th^\bot}\We^K$ and $x'c'_1 c'_2\in \mb{}^{\Th\cup \Th^\bot}\We$. Set $\Xi':=\Th\cup K^\infty \cup \red{x'}$. Then
\begin{eqnarray*}
   \ve{R(\Th)} z\We_K = a'_1\We_{\Xi'\cup ((\Xi')^\bot\cap x' c'_1 K^0)}.
\end{eqnarray*}
By our choices of the decompositions of $z$, it follows 
\begin{eqnarray}\label{deceq}
   a_1= a'_1\quad\mb{ and }\quad a_2= a'_2 \quad\mb{ and }\quad xc_1c_2= x' c'_1 c'_2. 
\end{eqnarray}
Therefore, we only have to show 
\begin{eqnarray*}
  \Xi\cup (\Xi^\bot \cap xc_1 J^0)\supseteq \Xi'\cup ((\Xi')^\bot\cap x' c'_1 K^0).
\end{eqnarray*}

We first show $\Xi\supseteq \Xi'$. Decompose $c_2=c_{21}c_{22}$ with $c_{21}\in\We_{J^\infty}$ and $c_{22}\in\We_{J^0}$. By (\ref{KJc}) we get 
$(K^\infty)^\bot\supseteq (J^\infty)^\bot\supseteq J^0$. With (\ref{deceq}) we find
\begin{eqnarray*}
   &&   \We_{\Th^\bot}xc_{21}\We_{K^\infty \cup (K^\infty)^\bot}= \We_{\Th^\bot}xc_1c_2 \We_{K^\infty\cup (K^\infty)^\bot} \\
   && =   \We_{\Th^\bot}  x' c'_1 c'_2 \We_{K^\infty\cup (K^\infty)^\bot} = \We_{\Th^\bot} x' \We_{K^\infty \cup (K^\infty)^\bot}.
\end{eqnarray*} 
It follows $xc_{21}\geq x'$ because $x'$ is a minimal double coset representative. Therefore, $\red{xc_{21}}= \red{x}\cup\red{c_{21}}\supseteq \red{x'}$, and we find
\begin{eqnarray*}
  \Xi=\Th\cup J^\infty\cup \red{x} \supseteq \Th\cup K^\infty \cup \red{x'}= \Xi'.
\end{eqnarray*}

Use (\ref{KJc}) to decompose
\begin{eqnarray}
  (\Xi')^\bot\cap x' c'_1 K^0 = \left((\Xi')^\bot\cap x' c'_1 (K^0\cap J^\infty)\right)\cup \left((\Xi')^\bot\cap x' c'_1 (K^0\cap J^0)\right).
\end{eqnarray}
Next we show $(\Xi')^\bot\cap x' c'_1 (K^0\cap J^\infty) \subseteq\Xi$. By $K\subseteq J \subseteq J^\infty\cup (J^\infty)^\bot$, equations (\ref{deceq}), and 
the formula for the stabilizer of a face we find
\begin{eqnarray}\label{eqf}
   && R(\Xi)=R(\Th)\cap x R(J^\infty )=R(\Th)\cap x c_1c_2 R(J^\infty ) = R(\Th)\cap x' c'_1  c'_2 R(J^\infty )\nonumber 
      \\ && =R(\Th)\cap x' c'_1 R(J^\infty).
\end{eqnarray}
Let $i\in (\Xi')^\bot\cap x' c'_1 (K^0\cap J^\infty)$. Then the simple reflection $\sigma_i$ fixes the face $x' c'_1 R(J^\infty)$ pointwise. Therefore, 
it also fixes $R(\Xi)=R(\Th)\cap x' c'_1 R(J^\infty)$ pointwise. The parabolic subgroup $\We_\Xi$ is the pointwise stabilizer of $R(\Xi)$. It follows $i\in\Xi$.

It remains to show $(\Xi')^\bot\cap x' c'_1 (K^0\cap J^0) \subseteq\Xi\cup (\Xi^\bot \cap x c_1 J^0)$. 
Let $i\in (\Xi')^\bot\cap x' c'_1 (K^0\cap J^0)$. 
The parabolic subgroup $\We_{\Th\cup\Th^\bot}$ fixes $R(\Th)$ as a whole, and $\Th\subseteq \Xi'$ implies $\Th^\bot\supseteq (\Xi')^\bot$. Therefore, the simple 
reflection $\sigma_i$ fixes $R(\Th)$ as a whole. Similarly, it fixes $x' c'_1 R(J^\infty )$ as a whole. By (\ref{eqf}), it fixes $R(\Xi)$ as a whole. 
The parabolic subgroup $\We_{\Xi\cup\Xi^\bot}$ is the stabilizer of $R(\Xi)$ as a whole. It follows $i\in\Xi\cup\Xi^\bot$. Now suppose that $i\in\Xi^\bot$.
Note that 
\begin{eqnarray}\label{JWJ}
   \Xi^\bot \cap x c_1 J^0 = \Xi^\bot \cap x c_1 \We_J J^0 .
\end{eqnarray}
Here the nontrivial inclusion ``$\supseteq$'' can be seen as follows. Let $k\in \Xi^\bot$ and $\sum_{j\in J^0} n_j \al_j\in \We_{J^0}J^0=\We_J J^0$ such that
\begin{eqnarray*}
  \al_k=xc_1\sum_{j\in J^0} n_j \al_j = \sum_{j\in J^0} n_j xc_1 \al_j.
\end{eqnarray*}
Here $xc_1 \al_j$, $j\in J^0$, are positive roots because of $xc_1\in \mb{}^{\Th\cup \Th^\bot}\We^J$. All $n_j$ are nonnegative integers or all 
$n_j$ are nonpositive integers. Therefore, only one of the integers $n_j$ is 1, the others are 0.

Recall that $K\subseteq J$, $ x' c'_1 c'_2= xc_1c_2$, and $\Xi\supseteq \Xi'$, which implies $\Xi^\bot \subseteq (\Xi')^\bot$. With (\ref{JWJ}) it follows
\begin{eqnarray*}
    && i\in \Xi^\bot\cap (\Xi')^\bot \cap  x' c'_1 (K^0\cap J^0)  =  \Xi^\bot \cap x' c'_1 (K^0\cap J^0) 
     \subseteq  \Xi^\bot \cap x' c'_1 \We_K (K^0\cap J^0) \\
     &&=  \Xi^\bot \cap x' c'_1 c'_2\We_K (K^0\cap J^0) 
     \subseteq \Xi^\bot \cap x c_1 c_2\We_J J^0=  \Xi^\bot \cap x c_1 \We_J J^0 = \Xi^\bot \cap x c_1 J^0.    
\end{eqnarray*}
\qed

\begin{rem} For $\sigma,\,\sigma'\in\We$ and $J,\,J'\subseteq I$ the equivalence
\begin{eqnarray*}
  \sigma\We_J\subseteq \sigma'\We_{J'}\;\iff\; \sigma\WeD_J\subseteq \sigma'\WeD_{J'} 
\end{eqnarray*} 
holds. Therefore, it is possible and looks natural to take as Coxeter complex the set
\begin{eqnarray*}
   \widetilde{\mathcal C}:=\Mklz{\sigma \WeD_J}{\sigma\in\We,\;J\subseteq I},
\end{eqnarray*}
partially ordered by the reverse inclusion. But multiplication on the left by the elements of $\WeD$ does not induce an action of $\WeD$ on 
$\widetilde{\mathcal C}$. 
The set $\widetilde{\mathcal C}$ is in general not stable by these operations. For example, let $A$ be of indefinite, but not strongly hyperbolic 
type. Then there exists a nontrivial special set $\emptyset\neq \Th\neq I$, and 
\begin{eqnarray*}
  e(R(\Th))\WeD_\emptyset =\{ e(R(\Th)) \} \neq \sigma \WeD_J \;\mb{ for all }\; \sigma\in\We,\;J\subseteq I.
\end{eqnarray*}
\end{rem}
%
%
%
%%%%%%%%%%%%%%%%%%%%%%%%%%%%%%%%%%%%%%%%%%%%%%%%%%%%%%%%%%%%%%%%%%%%%%%%%%%%%%%%%%%%%%%%%%%%%%%%%%%%%%%%%%%%%%%%%%%%%%%%%%%%%%%%%%%%%%%%
%
%
\section{The corresponding actions of the face monoid $\GD$ associated to $G$ on the building $\Omega$}
%
%
%%%%%%%%%%%%%%%%%%%%%%%%%%%%%%%%%%%%%%%%%%%%%%%%%%%%%%%%%%%%%%%%%%%%%%%%%%%%%%%%%%%%%%%%%%%%%%%%%%%%%%%%%%%%%%%%%%%%%%%%%%%%%%%%%%%%%%%%
% 
%

The Kac-Moody group $G$ acts naturally on its building $\Omega$. In this section we extend this action to actions of the face monoid $\GD$ on $\Omega$, which fit to the actions of the face monoid $\WeD$ on the Coxeter complex $\mathcal C$ given in the last section. To describe the stabilizers of our actions we need the analogues of the parabolic subgroups of $G$, the parabolic submonoids of $\GD$, which we treat first. 
\begin{thm}\label{AB1} For $J\subseteq I$ set
\begin{eqnarray*}
   \LD_J:= (T U_J) \WeD_J (T U_J)  \quad\mb{ and }\quad\PD_J := B \,\WeD_J B.
\end{eqnarray*}
Then $\LD_J$ and $\PD_J$ are submonoids of $\GD$ with unit groups 
\begin{eqnarray*}\label{Leviint}
   (\LD_J)^\times=\LD_J\cap G=L_J \quad\mb{ and }\quad (\PD_J)^\times=\PD_J\cap G=P_J.
\end{eqnarray*} 
The multiplication map $m:\LD_J\times U^J\to \PD_J$ is bijective. 
\end{thm}

We call $\PD_J$ the {\it standard parabolic submonoid} of $\GD$ of {\it type }$J$, and every conjugate $g\PD_J g^{-1}$, $g\in G$, a {\it parabolic submonoid} of $\GD$ of {\it type} $J$. We call $\PD_J=\LD_J U^J$ the {\it standard Levi decomposition} of $\PD_J$ with {\it Levi part} $\LD_J$ and {\it unipotent part} $U^J$.
\begin{rem}\label{AB2}
In the standard Levi decomposition $\PD_J=\LD_J U^J$ the order of the factors can not be changed: It holds $U^J\LD_J\subseteq \PD_J$, but in general this inclusion is proper.
\end{rem}
\begin{rem}\label{AB3} Equip $\GD$ with the Zariski topology induced by its coordinate ring, compare \cite{M1}. Then 
\begin{eqnarray*}
  \LD_J\subseteq \overline{L_J}\;\mb{ and }\;\PD_J\subseteq \overline{P_J},
\end{eqnarray*} 
but in general these inclusions are proper. For example, $\LD_\emptyset=\PD_\emptyset=B$ but 
\begin{eqnarray*}
\overline{B}=\dot{\bigcup_{\Th\; special}}\;\;\dot{\bigcup_{\sigma_1\in{\mathcal W}^\Th,\;\sigma_2\in\mb{}^{\Th\cup\Th^\bot}{\mathcal W}\atop \sigma_1\leq\sigma_2^{-1} }} 
            B\sigma_1 e(R(\Th))\sigma_2 B
\end{eqnarray*}
by Theorem 9 of \cite{M3}.
\end{rem}
 \Proof By Section 2.9 of \cite{M1}, $\LD_J = B_J \WeD_J B_J =T \GD_J=\GD_J T$ where $\GD_J\cong \widehat{G(A_J)'}$. In particular, $\LD_J$ is a submonoid of $\GD$ with unit group $\LD_J\cap G=L_J$.

We next show $\PD_J=\LD_J U^J$. Trivially, 
\begin{eqnarray*}
   \PD_J = B\WeD_J B = U^J (T U_J) \WeD_J (T U_J) U^J\supseteq \LD_J U^J .
\end{eqnarray*}
To show the reverse inclusion note that $\PD_J=U_J U^J\WeD_J B$ and $U_J\subseteq \LD_J$. Therefore, it is sufficient to show $U^J T\WeD_J\subseteq \WeD_J B$.
Let $\sigma_1 \ve{R(\Th)} \sigma_2 \in\WeD_J$ be in normal form I, i.e., $\Th\subseteq J$ special, $\sigma_1\in(\We_J)^\Th=\We_J\cap \We^\Th$, and 
$\sigma_2\in\mb{}^{\Th\cup (\Th^\bot\cap J)}(\We_J)=\We_J\cap \mb{}^{\Th\cup\Th^\bot}\We$. 
Let $n_{\sigma_1}e(R(\Th)) n_{\sigma_2}\in\ND$ project onto $\sigma_1 \ve{R(\Th)} \sigma_2 \in\WeD=\ND/T$. Then for $u\in U^J$ we find
\begin{eqnarray*}
   && u n_{\sigma_1}e(R(\Th)) n_{\sigma_2} = n_{\sigma_1} \underbrace{(n_{\sigma_1})^{-1} u n_{\sigma_1}}_{\in U^J} e(R(\Th)) n_{\sigma_2}
     = n_{\sigma_1}   e(R(\Th)) \underbrace{ p_\Th((n_{\sigma_1})^{-1} u n_{\sigma_1})}_{\in U_{\Th^\bot}} n_{\sigma_2} \\
    && = n_{\sigma_1}   e(R(\Th))n_{\sigma_2}\underbrace{(n_{\sigma_2})^{-1} p_\Th((n_{\sigma_1})^{-1} u n_{\sigma_1}) {n_{\sigma_2}}}_{\in U}. 
\end{eqnarray*}
Here we used that $U^J$ is normal in $P_J$, Theorem \ref{FM1} (b), and $(\sigma_2)^{-1}\in\We^{\Th^\bot}$.

$\PD_J$ is a submonoid of $\GD$, because of $1\in \PD_J$ and
\begin{eqnarray*}
  \PD_J \PD_J = (\LD_J U^J)( \LD_J U^J)  = \LD_J ( U^J  \LD_J U^J) = \LD_J \PD_J = \LD_J \LD_J U^J =\LD_J U^J =\PD_J.
\end{eqnarray*}
Obviously, its unit group is $\PD_J\cap G= P_J$.

It remains to show that the map $m$ is injective. Suppose that
\begin{eqnarray}\label{injgl}
  u n_\tau e(R(\Th)) n_\sigma v x  =   u' n'_{\tau'} e(R(\Th')) n'_{\sigma'} v' x',
\end{eqnarray}
where $u,\,v,\, u',\, v'\in U_J$, $x,\, x'\in U^J$, and $n_\tau  e(R(\Th)) n_\sigma,\, n'_{\tau'} e(R(\Th')) n'_{\sigma'}\in
\ND_J=\WeD_J T$. We have to show $u n_\tau e(R(\Th)) n_\sigma v =   u' n'_{\tau'} e(R(\Th')) n'_{\sigma'} v'$ and $x= x'$. 
Since $U^J$ is a group it is sufficient to show $x= x'$.

By the Bruhat decomposition of $\GD$ from (\ref{injgl}) follows
\begin{eqnarray}\label{injwdgl}
   n_\tau e(R(\Th)) n_\sigma  = n'_{\tau'} e(R(\Th')) n'_{\sigma'}.
\end{eqnarray}
We may assume that $\tau \ve{R(\Th)}\sigma\in\WeD_J$ is in normal form I, i.e., $\Th\subseteq J$ special, $\tau\in(\We_J)^\Th=\We_J\cap \We^\Th$, 
and $\sigma\in \mb{}^{\Th\cup(\Th^\bot \cap J)}(\We_J)=\mb{}^{\Th\cup\Th^\bot}\We\cap\We_J$. 
We may also assume that $ u n_\tau e(R(\Th)) n_\sigma v$ and $ u' n_\tau e(R(\Th)) n_\sigma v'$ are in normal form, i.e., 
$u,\,u'\in (U_J)_\tau= U_J \cap \tau ((U_J)^\Th)^- \tau^{-1}$, and $v,\,v'\in U_J\cap \sigma^{-1} (U_J)^\Th \sigma$. 
(Use $\LD_J=T \GD_J=\GD_J T$ and $\GD_J\cong \widehat{\GD(A_J)'}$ and Theorem \ref{FM6}.) By (\ref{injgl}) and (\ref{injwdgl}) we get
\begin{eqnarray}
 \underbrace{ {n_\tau}^{-1} u^{-1} u' n_\tau}_{\in ((U_J)^\Th)^-} e(R(\Th)) = 
  e(R(\Th))n_\sigma v x (x')^{-1} (v')^{-1} {n_\sigma}^{-1} \hspace{5em}\nonumber\\
   =e(R(\Th))\underbrace{ \underbrace{n_\sigma v {n_\sigma}^{-1}}_{\in (U_J)^\Th} \underbrace{n_\sigma x (x')^{-1} {n_\sigma}^{-1}}_{\in U^J}  
  \underbrace{n_\sigma (v')^{-1} {n_\sigma}^{-1}}_{\in(U_J)^\Th}}_{\in U^\Th}.\label{injugl}
\end{eqnarray} 
To cut short the notation set $a:={n_\tau}^{-1} u^{-1} u' n_\tau$ and $b:=n_\sigma v x (x')^{-1} (v')^{-1} {n_\sigma}^{-1}$. From (\ref{injugl}) follows
\begin{eqnarray*}
  a\in P_{\Th\cup\Th^\bot}\cap ((U_J)^\Th)^- \subseteq (U_{\Th\cup\Th^\bot})^-
\end{eqnarray*} 
and 
\begin{eqnarray*}
  b\in P_{\Th\cup\Th^\bot}^-\cap U^\Th= U_{\Th\cup\Th^\bot}\cap U^\Th=U_{\Th^\bot}.
\end{eqnarray*}
(Here we used the fact that $P_K\cap U^-=U_K^-$ and $P_K^-\cap U =U_K$, $K\subseteq I$. This can be checked easily. Use the Birkhoff decompositions of $P_K^-$, $P_K$, $G$, 
and $U\cap U^-=\{1\}$.)
Furthermore, $U^-\ni p_\Th(a)=p_\Th^-(b)=b\in U^+$. It follows $b=1$, which gives
\begin{eqnarray*}
    U^J\ni n_\sigma x (x')^{-1}{n_\sigma}^{-1} = n_\sigma v^{-1} {n_\sigma}^{-1}\,  1 \, n_\sigma  v' {n_\sigma}^{-1} \in (U_J)^\Th.
\end{eqnarray*}  
Since $U^J\cap (U_J)^\Th= \{1\}$ this shows $x=x'$.
\qed 
\begin{prop}\label{AB4} Let $J,\,K\subseteq I$. Then:
\begin{itemize}
\item[(a)] $\PD_J\subseteq \PD_K$ if and only if $J\subseteq K$.
\item[(b)] $\PD_J\cap\PD_K=\PD_{J\cap K}$.
\item[(c)] $\PD_J$ and $\PD_K$ generate $\PD_{J\cup K}$ as a monoid.
\end{itemize}
\end{prop}
 \Proof Part (a) and (b) follow immediately from part (a) and (b) of Proposition \ref{AC1}, use the Bruhat decomposition of $\GD$. Part (c) follows from part (c) of Proposition \ref{AC1}.\qed 

We take as building the set
\begin{eqnarray*}
    \Omega := \Mklz{g P_J}{g\in G,\;J\subseteq I}
\end{eqnarray*}
partially ordered by the reverse inclusion, i.e., for $g,\,g'\in G$ and $J,\,J'\subseteq I$, 
\begin{eqnarray*}
    g P_J\leq g' P_{J'} \;:\iff\; g P_J\supseteq g' P_{J'}. 
\end{eqnarray*}
The Kac-Moody group $G$ acts order preservingly on $\Omega$ by multiplication from the left. 

Let ${\mathcal A}:=\Mklz{n P_J}{n\in N,\;J\subseteq I}$ be the standard apartment of $\Omega$. Denote by $\pi:\ND\to\WeD$ the canonical projection of $\ND$ onto $\WeD=\ND/T$. As an ordered set the standard appartement $\mathcal A$ identifies with the Coxeter complex $\mathcal C$ by the map $\rho:{\mathcal A}\to{\mathcal C}$ defined by $\rho( n P_J ):=\pi(n)\We_J$, $n\in N$ and $J\subseteq I$. Set $\psi:=\rho^{-1}:{\mathcal C}\to {\mathcal A}$.
\begin{defn}\label{AB5} Let $\bullet$ be an action of the face monoid $\WeD$ on the Coxeter complex $\mathcal C$, which extends the action of the Weyl group $\We$ on $\mathcal C$. An action $\bullet$ of the face monoid $\GD$ on the building $\Omega$, which extends the action of the Kac-Moody group $G$ on $\Omega$, is called compatible if the following equivalent conditions hold:
\begin{itemize}
\item[(i)] It holds $\ND_\bullet {\mathcal A} \subseteq {\mathcal A}$ and for every $\hat{n}\in\ND$ the diagram
\begin{eqnarray*}
\begin{array}{rcl}
 {\mathcal A} &\stackrel{\rho}{\longrightarrow} &{\mathcal C}\\
  \hat{n}_\bullet \downarrow     &    & \downarrow\pi(\hat{n})_\bullet \\
 {\mathcal A} &\stackrel{\rho}{\longrightarrow} & {\mathcal C}
\end{array}
\end{eqnarray*}
commutes.
\item[(ii)] For all special sets $\Th\subseteq I$ it holds $e(R(\Th))_\bullet {\mathcal A} \subseteq {\mathcal A}$ and the diagram 
\begin{eqnarray*}
\begin{array}{rcl}
 {\mathcal A} &\stackrel{\rho}{\longrightarrow} &{\mathcal C}\\
  e(R(\Th))_\bullet \downarrow     &    & \downarrow\ve{R(\Th)}_\bullet \\
 {\mathcal A}&\stackrel{\rho}{\longrightarrow} & {\mathcal C}
\end{array}
\end{eqnarray*}
commutes.
\end{itemize}
\end{defn}
\begin{rem}\label{AB6} For actions as in the definition the action of $\GD$ on $\Omega$ can be recovered from the action of $\WeD$ on $\mathcal C$ as follows:

Let $\hat{g}\in\GD$ and $hP_J\in\Omega$. Write $\hat{g}$ in the form $\hat{g}=g_1 e(R(\Th)) g_2$ with 
$g_1,\,g_2\in G$ and $\Th$ special. Decompose $g_2 h= a n_\sigma b$ with $a\in P_{\Th\cup\Th^\bot}^-$, $n_\sigma\in N$ projecting to $\sigma\in\We$, and $b\in P_J$.
Then
\begin{eqnarray*}
   && \hat{g}_\bullet h P_J = g_1 e(R(\Th))g_2\mb{}_\bullet hP_J = g_1 e(R(\Th))a_\bullet n_\sigma b P_J\\
   && = g_1 p_\Th^-(a) e(R(\Th))_\bullet n_\sigma P_J = g_1 p_\Th^-(a) \psi(\ve{R(\Th)}_\bullet \sigma \We_J ). 
\end{eqnarray*}
In particular, there exists at most one compatible action of the face monoid $\GD$ on the building $\Omega$, which extends the action of the Kac-Moody group $G$ on $\Omega$.
\end{rem}

We now fix an action $\bullet$ of the face monoid $\WeD$ on the Coxeter complex $\mathcal C$ which extends the action of the Weyl group $\We$ on $\mathcal C$. We determine necessary and sufficient conditions such that the formula of the last remark can be taken as definition of a compatible action of the face monoid $\GD$ on the building $\Omega$.
As a first step we show:
\begin{thm}\label{AB7} Let $\Th\subseteq I$ be a special set. The following statements are equivalent:
\begin{itemize}
\item[(i)] We get a well defined map $e(R(\Th))_\bullet:\Omega \to\Omega$ in the following way: Let $gP_J\in\Omega$. Decompose $g$ in the form $g=an_\sigma b$ 
with $a\in P_{\Th\cup\Th^\bot}^-$, $n_\sigma\in N$ projecting to $\sigma\in\We$, and $b\in P_J$. Define
\begin{eqnarray*}
  e(R(\Th))_\bullet gP_J:= p_\Th^-(a)\psi(\ve{R(\Th)}_\bullet\sigma \We_J).
\end{eqnarray*}
\item[(ii)] $p_\Th^-(P_{\Th\cup\Th^\bot}^-\cap\sigma P_J\sigma^{-1})\psi(\ve{R(\Th)}_\bullet\sigma\We_J)=\psi(\ve{R(\Th)}_\bullet\sigma\We_J)$ for all 
$\sigma\in\mb{}^{\Th\cup\Th^\bot}\We^J$, $J\subseteq I$.
\end{itemize}
If this holds then the map $e(R(\Th))_\bullet:\Omega\to\Omega $ is order preserving if and only if the map $\ve{R(\Th)}_\bullet:{\mathcal C}\to{\mathcal C}$ 
is order preserving. 
\end{thm}
\Proof Let $gP_J=g'P_{J'}$. Then $J'=J$ and $g^{-1}g'\in P_J$. Decompose $g=an_\sigma b$ and $g'=a'n'_{\sigma'}b'$ as in the Theorem.
Then $ e(R(\Th))_\bullet$ applied to  $gP_J=g'P_{J'}$ is well defined if and only if 
\begin{eqnarray*}
p_\Th^-(a)\psi(\ve{R(\Th)}_\bullet\sigma\We_J)=p_\Th^-(a')\psi(\ve{R(\Th)}_\bullet\sigma'\We_J).
\end{eqnarray*}
Now decompose $\sigma = \sigma_1\sigma_2 \sigma_3\sigma_4$ and 
$\sigma' = \sigma'_1\sigma'_2 \sigma'_3\sigma'_4$ with $\sigma_1,\,\sigma'_1\in \We_{\Th^\bot}$, $\sigma_2,\,\sigma'_2\in\We_\Th$, 
$\sigma_3,\,\sigma'_3\in \mb{}^{\Th\cup\Th^\bot}\We^J$, and $\sigma_4,\,\sigma'_4\in \We_J$. Then
\begin{eqnarray*}
  &&   p_\Th^-(a)\psi(\ve{R(\Th)}_\bullet\sigma\We_J)=p_\Th^-(a)\psi(\ve{R(\Th)}_\bullet\sigma_1\sigma_2\sigma_3\sigma_4\We_J)\\
  && = p_\Th^-(a)\psi(\ve{R(\Th)}\sigma_1\sigma_2\mb{}_\bullet\sigma_3\sigma_4\We_J) = p_\Th^-(a)\psi(\sigma_1\ve{R(\Th)}\mb{}_\bullet\sigma_3\We_J)\\ 
  && = p_\Th^-(a)\widetilde{\sigma_1}\psi(\ve{R(\Th)}\mb{}_\bullet\sigma_3 \We_J) = p_\Th^-(a\widetilde{\sigma_1})\psi(\ve{R(\Th)}\mb{}_\bullet\sigma_3 \We_J)\\
  && = p_\Th^-(a\widetilde{\sigma_1}\widetilde{\sigma_2})\psi(\ve{R(\Th)}\mb{}_\bullet\sigma_3 \We_J).
\end{eqnarray*}
Therefore, $e(R(\Th))_\bullet$ applied to  $gP_J=g'P_{J'}$ is well defined if and only if 
\begin{eqnarray*}
         p_\Th^-(a\widetilde{\sigma_1}\widetilde{\sigma_2})\psi(\ve{R(\Th)}\mb{}_\bullet\sigma_3 \We_J)  
         =   p_\Th^-(a'\widetilde{\sigma'_1}\widetilde{\sigma'_2})\psi(\ve{R(\Th)}\mb{}_\bullet\sigma'_3 \We_J).   
\end{eqnarray*}
Now $gP_J=a n_\sigma P_J=(a \widetilde{\sigma_1}\widetilde{\sigma_2})\widetilde{\sigma_3}P_J$, similarly 
$g'P_J= (a' \widetilde{\sigma'_1}\widetilde{\sigma'_2})\widetilde{\sigma'_3}P_J$,  where 
$a \widetilde{\sigma_1}\widetilde{\sigma_2},\,a' \widetilde{\sigma'_1}\widetilde{\sigma'_2}\in P_{\Th\cup\Th^\bot}^-$. Therefore, from $gP_J=g'P_J$ we get
$\sigma'_3=\sigma_3$ and $(a\widetilde{\sigma_1}\widetilde{\sigma_2})^{-1}a'\widetilde{\sigma'_1}\widetilde{\sigma'_2}\in \sigma_3 P_J(\sigma_3)^{-1}$.
It follows that $e(R(\Th))_\bullet$ applied to  $gP_J=g'P_{J'}$ is well defined if and only if 
\begin{eqnarray}\label{fwell}
         p_\Th^-((a\widetilde{\sigma_1}\widetilde{\sigma_2})^{-1}a'\widetilde{\sigma'_1}\widetilde{\sigma'_2})\psi(\ve{R(\Th)}\mb{}_\bullet\sigma_3 \We_J)  
         =  \psi(\ve{R(\Th)}\mb{}_\bullet\sigma_3 \We_J).   
\end{eqnarray}
Furthermore, here $(a\widetilde{\sigma_1}\widetilde{\sigma_2})^{-1}a'\widetilde{\sigma'_1}\widetilde{\sigma'_2}\in P_{\Th\cup\Th^\bot}^-\cap\sigma_3 P_J(\sigma_3)^{-1}$.

Depending on $g$ and $g'$, all the elements of the intersection $P_{\Th\cup\Th^\bot}^-\cap\sigma_3 P_J(\sigma_3)^{-1}$ appear. This can be seen as follows. Let 
$\sigma_3\in \mb{}^{\Th\cup\Th^\bot}\We^J$ and $x\in P_{\Th\cup\Th^\bot}^-\cap\sigma_3 P_J(\sigma_3)^{-1}$. Then $\widetilde{\sigma_3}P_J=x\widetilde{\sigma_3}P_J$. Choose 
$a=1$, $a'=x$, $\sigma_1=\sigma'_1=1$,  $\sigma_2=\sigma'_2=1$, and $\sigma_4=\sigma'_4=1$. Then (\ref{fwell}) reduces to    
\begin{eqnarray*}
         p_\Th^-(x)\psi(\ve{R(\Th)}\mb{}_\bullet\sigma_3 \We_J)  =  \psi(\ve{R(\Th)}\mb{}_\bullet\sigma_3 \We_J).   
\end{eqnarray*}

This proves the equivalence of (i) and (ii). Now suppose that (i), (ii) hold. It remains to show that $e(R(\Th))_\bullet$ is order preserving if and only if $\ve{R(\Th)}_\bullet$ is order preserving. 

It holds $\ve{R(\Th)}_\bullet =\rho\circ e(R(\Th))_\bullet \circ \psi$ by the definition of $e(R(\Th))_\bullet$. The maps $\psi$ and $\rho$ are order 
preserving. Therefore, if $e(R(\Th))_\bullet$ is order preserving, then also $\ve{R(\Th)}_\bullet$ is order preserving.

Now suppose that $\ve{R(\Th)}_\bullet$ is order preserving. Let $g'P_J,\,gP_K\in\Omega$ with $g'P_J\leq g P_K$. This is equivalent to 
$g'P_J\supseteq gP_K$, which in turn is equivalent to $J\supseteq K$ and $(g')^{-1}g\in P_J$. In particular, $g'P_J=g P_J$. Now decompose $g=a n_\sigma b$ with
$a\in P_{\Th\cup\Th^\bot}$ $n_\sigma\in N$ projecting to $\sigma\in\We$, and $b\in P_K\subseteq P_J$. Then
\begin{eqnarray*}
  &&e(R(\Th))_\bullet g' P_J=e(R(\Th))_\bullet g P_J = p_\Th^-(a)\psi(\ve{R(\Th)}_\bullet\sigma\We_J),\\
  &&e(R(\Th))_\bullet g P_K = p_\Th^-(a)\psi(\ve{R(\Th)}_\bullet\sigma\We_K).
\end{eqnarray*}
It holds $\sigma \We_J\leq \sigma\We_K$. Furthermore,  $\ve{R(\Th)}_\bullet$ is order preserving, the map $\psi$ is order preserving, and $G$ acts order 
preservingly on $\Omega$. Therefore,
\begin{eqnarray*}
    e(R(\Th))_\bullet g' P_J= p_\Th^-(a)\psi(\ve{R(\Th)}_\bullet\sigma\We_J)
    \leq  p_\Th^-(a)\psi(\ve{R(\Th)}_\bullet\sigma\We_K)  = e(R(\Th))_\bullet g P_K.
\end{eqnarray*}
\qed

Next we want to determine the groups $p_\Th^-(P_{\Th\cup\Th^\bot}^-\cap\sigma P_J\sigma^{-1})$, $\sigma\in\mb{}^{\Th\cup\Th^\bot}\We^J$, which appear in (ii) of the last theorem, explicitely.
\begin{thm}\label{AB8} Let $J,\,K\subseteq I$. Let $\sigma\in \mb{}^K\We^J$. Then
\begin{eqnarray*}
   P_K^-\cap\sigma P_J\sigma^{-1} = (\sigma U_J^-\sigma^{-1}) (\sigma U^J\sigma^{-1}\cap (U^K)^-) N_{K\cap \sigma J} U_K.
\end{eqnarray*}
\end{thm}
 \Proof Since $\sigma\in \mb{}^K\We^J$ the inclusion ``$\supseteq$'' holds. It remains to show the reverse inclusion. By the Birkhoff decompositions of $P_K^-$, $P_J$ it holds
\begin{eqnarray*}
   P_K^-\ti{\sigma}\cap \ti{\sigma}P_J = U^-\We_K T U_K \ti{\sigma}\cap \ti{\sigma}U_J^-\We_J T U = U^-\We_K \sigma T (\sigma^{-1}U_K \sigma) \cap (\sigma U_J^-\sigma^{-1})\sigma\We_J T U. 
\end{eqnarray*}
Because of $\sigma\in \mb{}^K\We^J$ we have $\sigma^{-1}U_K\sigma\subseteq U$ and $\sigma U_J^-\sigma^{-1}\subseteq U^-$. By the Birkhoff decomposition of $G$ we get
\begin{eqnarray}\label{VPS}
       P_K^-\ti{\sigma}\cap \ti{\sigma}P_J= \bigcup_{\tau\in {\mathcal V}} 
              \left( U^- \tau T (\sigma^{-1}U_K \sigma) \cap (\sigma U_J^-\sigma^{-1})\tau T U  \right), 
\end{eqnarray}
where ${\mathcal V}:=\We_K\sigma\cap\sigma\We_J$. By Proposition 2.7.5. of \cite{Ca}, which is valid also for infinite Coxeter groups, it holds
\begin{eqnarray}\label{VWS}
   {\mathcal V} =\We_K\sigma\cap\sigma\We_J = (\We_K\cap \sigma \We_J\sigma^{-1})\sigma= \We_{K\cap\sigma J}\sigma . 
\end{eqnarray}

Let $x\in  P_K^-\ti{\sigma}\cap \ti{\sigma}P_J$. Then by (\ref{VPS}) there exist $u^-\in U^-$, $n_\tau\in N$ projecting to $\tau\in {\mathcal V}$, $u_K\in U_K$ such that $x=u^- n_\tau \ti{\sigma}^{-1}u_K\ti{\sigma}$. Since $u^-n_\tau \ti{\sigma}^{-1}u_K\ti{\sigma}\in\ti{\sigma}P_J$ and $\ti{\sigma}^{-1}u_K\ti{\sigma}\in U$ we get $u^- n_\tau\in \ti{\sigma}P_J$. Because of $\tau\in{\mathcal V}\subseteq \sigma\We_J$ we find $\ti{\sigma}^{-1}u^-\ti{\sigma}\in P_J$. It follows $u^-\in \sigma P_J\sigma^{-1}\cap U^-$. Therefore, we have found
\begin{eqnarray}\label{SPNP1}
   P_K^-\ti{\sigma}\cap \ti{\sigma}P_J\subseteq (\sigma P_J\sigma^{-1}\cap U^-) {\mathcal V}T \sigma^{-1} U_K \sigma.
\end{eqnarray}
Now we show 
\begin{eqnarray}\label{SPNP2}
   \sigma P_J\sigma^{-1}\cap U^-=\sigma U^-_J\sigma^{-1} (\sigma U\sigma^{-1}\cap U^-).
\end{eqnarray} 
By the Birkhoff decomposition for $P_J$ it holds
\begin{eqnarray*}
  \ti{\sigma}P_J\cap U^-\ti{\sigma}=(\sigma U_J^-\sigma^{-1})\sigma\We_J T U \cap U^-\ti{\sigma}.
\end{eqnarray*}
Here $(\sigma U_J^-\sigma^{-1})\subseteq U^-$ because of $\sigma\in \mb{}^K\We^J$. By the Birkhoff decomposition for $G$ we find
\begin{eqnarray*}
  \ti{\sigma}P_J\cap U^-\ti{\sigma}=(\sigma U_J^-\sigma^{-1})\ti{\sigma} U \cap U^-\ti{\sigma}.
\end{eqnarray*}
Multiplying on the right by $\ti{\sigma}^{-1}$, and using $\sigma U_J^-\sigma^{-1}\subseteq U^-$ we find 
\begin{eqnarray*}
  \sigma P_J\sigma^{-1}\cap U^-=(\sigma U_J^-\sigma^{-1})\sigma U \sigma^{-1}\cap U^- = (\sigma U_J^-\sigma^{-1})(\sigma U \sigma^{-1}\cap U^-) .
\end{eqnarray*}
Next we show 
\begin{eqnarray}\label{SPNP3}
        \sigma U\sigma^{-1}\cap U^-=\sigma U^J\sigma^{-1}\cap (U^K)^-.
\end{eqnarray}
The inclusion ``$\supseteq$'' is obvious, we only have to show the reverse inclusion. It holds
\begin{eqnarray*}
  \sigma U\sigma^{-1}\cap U^-=\prod_{\al\in\sigma\Delta^+_{re}\cap\Delta^-_{re}}U_\al.
\end{eqnarray*}
Let $\al\in\sigma\Delta^+_{re}\cap\Delta^-_{re}$. Suppose that $\al\in (\Delta^-_K)_{re}$. Since $\sigma\in \mb{}^K\We^J$ it follows $\sigma^{-1}\al\in\Delta^-_{re}$, which 
contradicts $\al\in\sigma\Delta^+_{re}$. Suppose that $\al\in \sigma (\Delta^+_J)_{re}$. Since $\sigma\in \mb{}^K\We^J$ it follows 
$\al\in\Delta^+_{re}$, which contradicts $\al\in\Delta^-_{re}$. Therefore
\begin{eqnarray*}
  \sigma U\sigma^{-1}\cap U^-\subseteq \prod_{\al\in\sigma(\Delta^+_{re}\setminus (\Delta^+_{J})_{re})\cap (\Delta^-_{re}\setminus (\Delta^-_K)_{re})}U_\al
\subseteq \sigma U^J\sigma^{-1}\cap (U^K)^-.
\end{eqnarray*}

Now multiply (\ref{SPNP1}) with $\ti{\sigma}^{-1}$ from the right. With (\ref{SPNP2}), (\ref{SPNP3}), and (\ref{VWS}) follows the inclusion ``$\subseteq$'' of the theorem.
\qed 
\begin{cor}\label{AB9} Let $\Th\subseteq I$ be special, let $J\subseteq I$ and $\sigma\in \mb{}^{\Th\cup\Th^\bot}\We^J$. Then
\begin{eqnarray}\label{CorP}
    p_\Th^-(P_{\Th\cup\Th^\bot}^-\cap\sigma P_J\sigma^{-1}) = U_{\Th^\bot\cap\sigma J}^-\widetilde{\We}_{\Th^\bot\cap \sigma J}T^\Th U_{\Th^\bot}.
\end{eqnarray}
\end{cor}
\begin{rem}\label{AB10} The image of $p_\Th^-$ is $G_{\Th^\bot}T^{\Th\cup\Th^\bot}$, which is a group with BN-pair 
$G_{\Th^\bot}T^{\Th\cup\Th^\bot}\cap B=U_{\Th^\bot}T^\Th$ and 
$G_{\Th^\bot}T^{\Th\cup\Th^\bot}\cap N=\widetilde{\We}_{\Th^\bot}T^\Th$. The group (\ref{CorP}) of the corollary is the standard parabolic subgroup corresponding to 
$\Th^\bot\cap\sigma J$.
\end{rem}
 \Proof By the last Theorem and our definition of $p_\Th^-$ we get
\begin{eqnarray}\label{SP1}
   \lefteqn{p_\Th^-(P_{\Th\cup\Th^\bot}^- \cap \sigma P_J\sigma^{-1})}&&\nonumber\\
      &=& p_\Th^-(\sigma U_J^-\sigma^{-1}) p_\Th^-(\sigma U^J\sigma^{-1}\cap (U^{\Th\cup\Th^\bot})^-) 
   p_\Th^-(\widetilde{\We}_{(\Th\cup\Th^\bot)\cap \sigma J})p_\Th^-(T) p_\Th^-(U_{\Th\cup\Th^\bot})\nonumber\\
 &=& p_\Th^-(\sigma U_J^-\sigma^{-1})\, 1 \,
  \widetilde{\We}_{\Th^\bot\cap \sigma J} T^\Th U_{\Th^\bot}.
\end{eqnarray} 
Since $\sigma\in\mb{}^{\Th\cup\Th^\bot}\We^J$ we get $\sigma U_J^-\sigma^{-1}\subseteq U^-$. By the definition of $p_\Th^-$ it follows $p_\Th^-(\sigma U_J^-\sigma^{-1})\subseteq U_{\Th^\bot}^-$. 

Now the group $p_\Th^-(P_{\Th\cup\Th^\bot}^- \cap \sigma P_J\sigma^{-1})$ contains $U_{\Th^\bot} T^\Th$. Therefore, it has to be a standard parabolic of $G_{\Th^\bot}T^{\Th\cup\Th^\bot}$ of the form
\begin{eqnarray}\label{SP2}
  U_K^- \widetilde{\We}_K T^\Th U_{\Th^\bot}
\end{eqnarray} 
for some $K\subseteq \Th^\bot$. If we now use the Birkhoff decomposition of $G$ to compare (\ref{SP1}) with (\ref{SP2}) we get $K=\Th^\bot\cap \sigma J$. 
\qed 

\begin{thm}\label{AB11} Let $\Th\subseteq I$ be a special set. Suppose that 
\begin{eqnarray*}
 U_{\Th^\bot\cap\sigma J}^-\widetilde{\We}_{\Th^\bot\cap \sigma J}T^\Th U_{\Th^\bot}\psi(\ve{R(\Th)}_\bullet\sigma\We_J)=\psi(\ve{R(\Th)}_\bullet\sigma\We_J)
\end{eqnarray*} 
holds for all $\sigma\in\mb{}^{\Th\cup\Th^\bot}\We^J$, $J\subseteq I$. Then the following statements are equivalent:
\begin{itemize}
\item[(i)] For every $\hat{g}\in G e(R(\Th))G$ we get a well defined map $\hat{g}:\Omega \to\Omega$ in the following way: Let $cP_J\in\Omega$. 
Write $\hat{g}$ in the form $\hat{g}=x e(R(\Th))y$ with $x,y\in G$. Define
\begin{eqnarray*}
  \hat{g}_\bullet cP_J:= x(e(R(\Th)_\bullet ycP_J).
\end{eqnarray*}
\item[(ii)] $G_\Th U^{\Th\cup\Th^\bot} \psi(\ve{R(\Th)}_\bullet\sigma\We_J)=\psi(\ve{R(\Th)}_\bullet\sigma\We_J)$ for all 
$\sigma\in\mb{}^{\Th\cup\Th^\bot}\We^J$, $J\subseteq I$.
\end{itemize}
\end{thm}
\Proof Suppose that (i) holds. Then we find for all $a\in G_\Th U^{\Th\cup\Th^\bot}=\Mklz{g\in G}{ge(R(\Th))=e(R(\Th))}$, for all $\sigma\in\mb{}^{\Th\cup\Th^\bot}\We^J$, and for all $J\subseteq I$:
\begin{eqnarray*}
    &&  a \psi(\ve{R(\Th)}_\bullet\sigma\We_J)= a (e(R(\Th))_\bullet\widetilde{\sigma}P_J ) = (a e(R(\Th)))_\bullet\widetilde{\sigma}P_J \\
    && = e(R(\Th))_\bullet\widetilde{\sigma}P_J = \psi(\ve{R(\Th)}_\bullet\sigma\We_J).
\end{eqnarray*}

Now suppose that (ii) holds. Let $x,\,x',\,y,\,y'\in G$ such that $\hat{g}=xe(R(\Th))y=x'e(R(\Th))y'$. We have to show that
\begin{eqnarray*}
  x (e(R(\Th))_\bullet ycP_J )=x'(e(R(\Th))_\bullet y'cP_J )
\end{eqnarray*}
for all $cP_J\in\Omega$. Actually, it is sufficient to show that $xe(R(\Th))=e(R(\Th))y$, $x,\,y\in G$, implies
\begin{eqnarray*}
  x (e(R(\Th))_\bullet cP_J)=e(R(\Th))_\bullet ycP_J
\end{eqnarray*}
for all $cP_J\in\Omega$. 

Fix $cP_J\in\Omega$. Decompose $c=a\widetilde{\sigma}b$ with $a\in P_{\Th\cup\Th^\bot}^-$, $\sigma\in\mb{}^{\Th\cup\Th^\bot}\We^J$, $b\in P_J$. The equality 
$xe(R(\Th))=e(R(\Th))y$ is equivalent to $x\in P_{\Th\cup\Th^\bot}$, $y\in P_{\Th\cup\Th^\bot}^-$  and $p_\Th(x)=p_\Th^-(y)$. It follows
\begin{eqnarray*}
  && e(R(\Th))_\bullet ycP_J = p_\Th^-(ya)\psi(\ve{R(\Th)}_\bullet\sigma \We_J)\\ 
  && =p_\Th^-(y)p_\Th^-(a)\psi(\ve{R(\Th)}_\bullet\sigma \We_J) = p_\Th(x)p_\Th^-(a)\psi(\ve{R(\Th)}_\bullet\sigma \We_J) \\ 
  && = x p_\Th^-(a)(p_\Th^-(a))^{-1}x^{-1}p_\Th(x)p_\Th^-(a)\psi(\ve{R(\Th)}_\bullet\sigma \We_J).
\end{eqnarray*}
By the definitions of $p_\Th^\pm$ it holds $x^{-1}p_\Th(x)\in G_\Th\ltimes U^{\Th\cup\Th^\bot}$ and also  $(p_\Th^-(a))^{-1}x^{-1}p_\Th(x)p_\Th^-(a)\in G_\Th\ltimes U^{\Th\cup\Th^\bot}$. Since (ii) is valid it follows
\begin{eqnarray*}
   e(R(\Th))_\bullet ycP_J= x p_\Th^-(a)\psi(\ve{R(\Th)}_\bullet\sigma \We_J)=x(e(R(\Th))_\bullet cP_J).
\end{eqnarray*}
\qed

We have fixed an action $\bullet$ of the face monoid $\WeD$ on the Coxeter complex $\mathcal C$, which extends the action of the Weyl group $\We$ on $\mathcal C$. The next beautiful theorem is the main theorem in this section. It characterizes in an unexpected easy way when there exists a compatible action of the face monoid $\GD$ on the building $\Omega$. It also describes this action.
\begin{thm}\label{AB12} (a) There exists a compatible action $\bullet$ of $\GD$ on $\Omega$ extending the action of $G$ on $\Omega$ if and only if $\ve{R(\Th)}_\bullet\sigma\We_J\in {\mathcal C}$ is a standard parabolic for all $\Th\subseteq I$ special, $J\subseteq I$, and $\sigma\in\mb{}^{\Th\cup\Th^\bot}\We^J$.

This action of $\GD$ on $\Omega$ is uniquely determined and can be obtained in the following way. 
Let $\hat{g}\in\GD$ and $hP_J\in\Omega$. Write $\hat{g}$ in the form $\hat{g}=g_1 e(R(\Th)) g_2$ with 
$g_1,\,g_2\in G$ and $\Th$ special. Decompose $g_2 h= a n_\sigma b$ with $a\in P_{\Th\cup\Th^\bot}^-$, $n_\sigma\in N$ projecting to $\sigma\in\We$, and $b\in P_J$.
Then
\begin{eqnarray}\label{acmap}
   \hat{g}_\bullet h P_J = g_1 p_\Th^-(a) \psi(\ve{R(\Th)}_\bullet \sigma \We_J ). 
\end{eqnarray}

(b) Let $J\subseteq I$. The stabilizer of $P_J\in\Omega$ of the action $\bullet$ of $\GD$ on $\Omega$ is obtained from the stabilizer of $\We_J\in{\mathcal C}$ of the action 
$\bullet$ of $\WeD$ on $\mathcal C$ by
\begin{eqnarray}\label{stabcom}
  \mb{Stab}_{\widehat{G}}(P_J)= B\mb{Stab}_{\widehat{\mathcal W}}(\We_J)B.
\end{eqnarray}
It is equivalent:
\begin{itemize}
\item[(i)] $\mb{Stab}_{\widehat{G}}(P_J)=\PD_J$.
\item[(ii)] $\mb{Stab}_{\widehat{\mathcal W}}(\We_J)=\WeD_J$.
\item[(iii)] $\mb{Stab}_{\widehat{\mathcal W}}(\We_J)\supseteq\WeD_J$.
\item[(iv)] $\ve{R(J^\infty)}_\bullet\We_J=\We_J$.
\end{itemize}

(c) The action $\bullet$ of $\GD$ on $\Omega$ is order preserving if and only if the action $\bullet$ of $\WeD$ on $\mathcal C$ is order preserving.
\end{thm}
\begin{rem} The equivalent conditions (i), (ii), (iii) and (iv) of (b) are not always satisfied for the actions which we consider. The easiest action of $\WeD$ on $\mathcal C$, which satisfies the required property stated in part (a) of the theorem, is obtained as follows: For $\sigma\We_J\in\mathcal C$ set
\begin{eqnarray*}
\begin{array}{lll}
\tau_\bullet \sigma\We_J := \tau\sigma\We_J &\mb{ for all }& \tau\in\We,\\
 \widehat{\tau}_\bullet \sigma\We_J := \We_I &\mb{ for all }& \widehat{\tau}\in\WeD\setminus\We.
\end{array}
\end{eqnarray*}
This action is even order preserving. But
\begin{eqnarray*}
  \mb{Stab}_{\widehat{\mathcal W}}(\We_J)=\left\{\begin{array}{lcl}
\We_J & \mb{if} & J\neq I\\
\WeD_I & \mb{if} &J=I
\end{array}\right. .
\end{eqnarray*}

Note also that the corresponding compatible action of $\GD$ on $\Omega$ is described as follows: For $gP_J\in \Omega$ it holds
\begin{eqnarray*}
\begin{array}{lll}
 h_\bullet g P_J= hg P_J &\mb{ if }& h\in G,\\
 \widehat{h}_\bullet gP_J= P_I &\mb{ if }& \widehat{h}\in\GD\setminus G.
\end{array}
\end{eqnarray*}
\end{rem}
\Proof ($\alpha$) We first check that the conditions (ii) of Theorem \ref{AB7} and \ref{AB11} hold for all $\Th\subseteq I$ special if and only if $\ve{R(\Th)}_\bullet\sigma\We_J$ is a standard parabolic subgroup of $\We$ for all $\Th\subseteq I$ special, $J\subseteq I$, and $\sigma\in\mb{}^{\Th\cup\Th^\bot}\We^J$.

Fix $\Th\subseteq I$ special, $J\subseteq I$, and $\sigma\in\mb{}^{\Th\cup\Th^\bot}\We^J$. The corresponding equations of the conditions (ii) of Theorem \ref{AB7} and Theorem \ref{AB11}, which have to be satisfied, are 
\begin{eqnarray*}
    U_{\Th^\bot\cap\sigma J}^-\widetilde{\We}_{\Th^\bot\cap \sigma J}T^\Th U_{\Th^\bot}\psi(\ve{R(\Th)}_\bullet\sigma\We_J)=\psi(\ve{R(\Th)}_\bullet\sigma\We_J),\\
    U_\Th^- \widetilde{\We_\Th}T_\Th U_\Th U^{\Th\cup\Th^\bot} \psi(\ve{R(\Th)}_\bullet\sigma\We_J)=\psi(\ve{R(\Th)}_\bullet\sigma\We_J).
\end{eqnarray*}
Because of $P_{\Th\cup(\Th^\bot\cap\sigma J)}=U_{\Th\cup(\Th^\bot\cap\sigma J)}^- N_{\Th\cup(\Th^\bot\cap\sigma J)}U$, these two equations are equivalent to
\begin{eqnarray}\label{ppsi}
    P_{\Th\cup(\Th^\bot\cap\sigma J)} \psi(\ve{R(\Th)}_\bullet\sigma\We_J)=\psi(\ve{R(\Th)}_\bullet\sigma\We_J).
\end{eqnarray}
Here $\psi(\ve{R(\Th)}_\bullet\sigma\We_J)$ is a coset of the form $gP_K$. Now $P_{\Th\cup(\Th^\bot\cap\sigma J)}gP_K=gP_K$ is equivalent to 
$g^{-1}P_{\Th\cup(\Th^\bot\cap\sigma J)}g\subseteq P_K$. By properties of the parabolic subgroups, compare Part 3.2.3 of \cite{Ti}, this is equivalent to 
$g\in P_K$ and $P_{\Th\cup(\Th^\bot\cap\sigma J)}\subseteq P_K$. Therefore, (\ref{ppsi}) is equivalent to $\psi(\ve{R(\Th)}_\bullet\sigma\We_J)$ is a standard 
parabolic subgroup of $G$ such that 
\begin{eqnarray*}
P_{\Th\cup(\Th^\bot\cap\sigma J)} \subseteq \psi(\ve{R(\Th)}_\bullet\sigma\We_J).
\end{eqnarray*} 
This in turn is equivalent to $\ve{R(\Th)}_\bullet\sigma\We_J$ is a standard parabolic subgroup of $\We$ with
\begin{eqnarray}\label{wpsi}
  \We_{\Th\cup(\Th^\bot\cap\sigma J)}\subseteq \ve{R(\Th)}_\bullet\sigma\We_J.
\end{eqnarray}

We now show that (\ref{wpsi}) is no proper restriction. It is satisfied for our action $\bullet$ of $\WeD$ on $\mathcal C$: Let $a\in \We_\Th$, $b\in \We_{\Th^\bot\cap\sigma J}=\We_{\Th^\bot}\cap \sigma\We_J\sigma^{-1}$. Then 
\begin{eqnarray*}
   a b\ve{R(\Th)}_\bullet\sigma\We_J=  a \ve{R(\Th)}b_\bullet  \sigma\We_J=   \ve{R(\Th)}_\bullet b \sigma\We_J= \ve{R(\Th)}_\bullet\sigma\We_J.
\end{eqnarray*}

($\beta$) From Theorem \ref{AB7} and Theorem \ref{AB11} now follows that for every $\hat{g}\in\GD$ we get by (\ref{acmap}) a well defined map $\hat{g}_\bullet:\Omega\to\Omega$. To show that we get an action of $\GD$, extending the action of $G$, it remains to show that for all $g_1,\,g_2,\, h_1,\, h_2\in G$, for all $\Th_1,\,\Th_2$ special and for all elements $zP_J\in\Omega$ it holds 
\begin{eqnarray*}
  (g_1 e(R(\Th_1)) g_2)_\bullet ((h_1 e(R(\Th_2)) h_2)_\bullet zP_J)=(g_1 e(R(\Th_1) g_2 h_1 e(R(\Th_2)) h_2)_\bullet zP_J. 
\end{eqnarray*}
Decompose $g_2 h_1=a \widetilde{\sigma} b$ with $a\in P_{\Th_1\cup\Th_1^\bot}^-$, $\sigma\in\mb{}^{\Th_1\cup\Th_1^\bot}\We^{\Th_2\cup\Th_2^\bot}$, $b\in P_{\Th_2\cup\Th_2^\bot}$. Then by (\ref{acmap}) and by Theorem \ref{FM1} (b) this equation is equivalent to 
\begin{eqnarray*}
    && g_1 p_{\Th_1}^-(a)e(R(\Th_1))_\bullet (\widetilde{\sigma} e(R(\Th_2))_\bullet p_{\Th_2}(b)h_2 zP_J )\\ 
    &&  \qquad\qquad\qquad  =   g_1 p_{\Th_1}^-(a)e(R(\Th_1\cup \Th_2\cup\red{\sigma }))_\bullet p_{\Th_2}(b)h_2 zP_J. 
\end{eqnarray*}
Therefore, it is sufficient to show
\begin{eqnarray}\label{sufact}
    e(R(\Th_1))_\bullet (\widetilde{\sigma} e(R(\Th_2))_\bullet zP_J ) 
                  =   e(R(\Th_1\cup \Th_2\cup\red{\sigma }))_\bullet zP_J
\end{eqnarray}
for all $\Th_1,\,\Th_2$ special, for all $\sigma\in \mb{}^{\Th_1\cup\Th_1^\bot}\We^{\Th_2\cup\Th_2^\bot}$, and for all $zP_J\in\Omega$.

Set $\Th:=\Th_1\cup\Th_2\cup\red{\sigma}$ for short. Decompose $z$ in the form $z=a \widetilde{x}b$ with $a\in P_{\Th^\bot}^-$, $x\in \We$ and $b\in P_J$. Then
\begin{eqnarray}\label{acgl1}
 e(R(\Th))_\bullet zP_J= p_\Th^-(a)\psi(\ve{R(\Th)}_\bullet x\We_J).
\end{eqnarray}
From $\Th_2\subseteq \Th$ follows $a\in P_{\Th^\bot}^-\subseteq P_{\Th_2^\bot}^-$. Therefore,
\begin{eqnarray}
    &&e(R(\Th_1))_\bullet(\ti{\sigma} e(R(\Th_2))_\bullet z P_J)= e(R(\Th_1))_\bullet \ti{\sigma} p_{\Th_2}^-(a)\psi(\ve{R(\Th_2)}_\bullet  x \We_J)\nonumber\\
    &&= e(R(\Th_1))_\bullet \ti{\sigma} p_{\Th_2}^-(a)\ti{\sigma}^{-1}\ti{\sigma}\psi(\ve{R(\Th_2)}_\bullet  x \We_J)\nonumber\\
    &&= e(R(\Th_1))_\bullet (\ti{\sigma} p_{\Th_2}^-(a)\ti{\sigma}^{-1})\psi(\sigma\ve{R(\Th_2)}_\bullet  x \We_J). \label{acd}
\end{eqnarray}
Suppose we can show
\begin{eqnarray}
 &&\ti{\sigma}p_{\Th_2}^-(a)\ti{\sigma}^{-1}\in P_{\Th_1^\bot}^-,\label{acdz1}\\
 &&p_{\Th_1}^-(\ti{\sigma}p_{\Th_2}^-(a)\ti{\sigma}^{-1}) e(R(\Th))=p_{\Th}^-(a) e(R(\Th)) \label{acdz2}.
\end{eqnarray}
Then (\ref{acd}) can be transformed by Theorem \ref{FM1} and Theorem \ref{FL4} in the following way:
\begin{eqnarray}
 &&   e(R(\Th_1))_\bullet (\ti{\sigma} p_{\Th_2}^-(a)\ti{\sigma}^{-1})\psi(\sigma\ve{R(\Th_2)}_\bullet  x \We_J)\nonumber\\
 && = p_{\Th_1}^-(\ti{\sigma}p_{\Th_2}^-(a)\ti{\sigma}^{-1})e(R(\Th_1))_\bullet \psi(\sigma\ve{R(\Th_2)}_\bullet x  \We_J)\nonumber\\
 && = p_{\Th_1}^-(\ti{\sigma}p_{\Th_2}^-(a)\ti{\sigma}^{-1})\psi(\ve{R(\Th_1)}_\bullet(\sigma\ve{R(\Th_2)}_\bullet  x \We_J))\nonumber\\
  && = p_{\Th_1}^-(\ti{\sigma}p_{\Th_2}^-(a)\ti{\sigma}^{-1})\psi(\ve{R(\Th)}_\bullet  x \We_J)\nonumber\\
  && = p_{\Th_1}^-(\ti{\sigma}p_{\Th_2}^-(a)\ti{\sigma}^{-1}) e(R(\Th))_\bullet\psi(x \We_J)\nonumber\\
  && = p_\Th^- (a) e(R(\Th))_\bullet \psi(x \We_J) = p_\Th^- (a) \psi(\ve{R(\Th)}_\bullet x \We_J).\label{acgl2}
\end{eqnarray}
The right side of (\ref{acgl1}) coincides with the right side of (\ref{acgl2}), i.e., equation (\ref{sufact}) is valid.

Now we show (\ref{acdz1}). Write $a\in P_{\Th^\bot}^-$ in the form $a= t \prod_{i=1,\,\ldots,\,m}u_{\beta_i}$ with $t\in T$, $u_{\beta_i}\in U_{\beta_i}$, $\beta_i\in \Delta^-_{re}\cup(\Delta^+_{\Th^\bot})_{re}$, $i=1,\,\ldots,\,m$, $m\in\N$.
The inclusion $\Th_2\subseteq \Th$ implies $\Th^\bot\subseteq (\Th_2)^\bot$. Therefore,
\begin{eqnarray}\label{zgl1}
  \ti{\sigma}p_{\Th_2}^-(a)\ti{\sigma}^{-1} =  
       \underbrace{\ti{\sigma}p_{\Th_2}^-(t)\ti{\sigma}^{-1}}_{\in T}\cdot
\prod_{i=1,\,\ldots,\,m\atop \beta_i\in (\Delta^-_{\Th_2^\bot})_{re}\cup(\Delta^+_{\Th^\bot})_{re}} \underbrace{\ti{\sigma}u_{\beta_i}\ti{\sigma}^{-1}}_{\in U_{\sigma\beta_i}} .
\end{eqnarray}
From $\red{\sigma}\subseteq \Th$ follows $\sigma\in \We_{\Th}$, from which in turn follows $\sigma (\Delta^+_{\Th^\bot})_{re}=(\Delta^+_{\Th^\bot})_{re}$.
Since $\sigma\in\mb{}^{\Th_1\cup\Th_1^\bot}\We^{\Th_2\cup\Th_2^\bot}$ we also get $\sigma (\Delta^-_{\Th_2^\bot})_{re}\subseteq \Delta^-_{re}$. Therefore, 
the expression on the right of (\ref{zgl1}) is contained in $P_{\Th^\bot}^-$. From the inclusion $\Th_1\subseteq \Th$ follows 
$P_{\Th^\bot}^-\subseteq P_{\Th_1^\bot}^-$.

Now we show (\ref{acdz2}). By Theorem \ref{FM1}, use (\ref{acdz1}), it holds: 
\begin{eqnarray*}
  && p_{\Th}^-(a)e(R(\Th))= e(R(\Th))a =e(R(\Th_1))\ti{\sigma} e(R(\Th_2))a \\
  && =p_{\Th_1}^-(\ti{\sigma}p_{\Th_2}^-(a)\ti{\sigma}^{-1})e(R(\Th_1))\ti{\sigma} e(R(\Th_2)) =  p_{\Th_1}^-(\ti{\sigma}p_{\Th_2}^-(a)\ti{\sigma}^{-1})e(R(\Th)).
\end{eqnarray*}

($\gamma$) Obviously, the action $\bullet$ of $\WeD$ is order preserving if and only if $\ve{R(\Th)}_\bullet:{\mathcal C}\to{\mathcal C}$ is order preserving for all $\Th\subseteq I$ special. Similarly, the action $\bullet$ of $\GD$ is order preserving if and only if $e(R(\Th))_\bullet:\Omega\to\Omega$ is order preserving for all $\Th\subseteq I$. Now part (c) of the theorem follows from Theorem \ref{AB7}.

($\delta$) Let $J\subseteq I$. Every element of $\hat{g}\in\GD$ can be written in the form $\hat{g}=b \hat{n}_{\hat{\sigma}} b'$ with $b,\,b'\in B$ and $\hat{n}_{\hat{\sigma}}\in \ND$ projecting to $\hat{\sigma}\in\WeD$. Now $P_J=\hat{g}_\bullet P_J=b n_{\hat{\sigma}} b'_\bullet P_J$ if and only if $(\hat{n}_{\hat{\sigma}})_\bullet P_J=P_J$ if and only if $\hat{\sigma}_\bullet \We_J=\We_J$. This shows formula (\ref{stabcom}) of part (b) of the theorem.

It remains to show the equivalence of (i), (ii), (iii) and (iv) of part (b). By formula (\ref{stabcom}) the statements (i) and (ii) are equivalent. By Proposition \ref{ACneu} the statements (iii) and (iv) are equivalent. Trivially, from (iii) follows (ii). Now suppose that the equivalent statements (iii) and (iv) hold. To show (ii) we have to show $\mb{Stab}_{\widehat{\mathcal W}}(\We_J)\subseteq\WeD_J$: Let $\hat{\sigma}\in \mb{Stab}_{\widehat{\mathcal W}}(\We_J)$. Write $\hat{\sigma}$ in normal form I, i.e., $\hat{\sigma}=\sigma_1\ve{R(\Th)} \sigma_2$ with $\sigma_1\in\We^\Th$, $\sigma_2\in\mb{}^{\Th\cup\Th^\bot}\We$, and $\Th\subseteq I$ special. Write $\sigma_2$ in the form $\sigma_2=xc_1c_2$ with $x\in \mb{}^{\Th\cup\Th^\bot}\We^{J^\infty\cup (J^\infty)^\bot}$, $c_1\in(\We_{(J^\infty)^\bot})^{J^0}$ and $c_2\in \We_J$ such that in addition $xc_1\in\mb{}^{\Th\cup\Th^\bot}\We^J$. Then $\We_J=\hat{\sigma}_\bullet\We_J$ is equivalent to
\begin{eqnarray}\label{stgg} 
   (\sigma_1)^{-1}\We_J=\ve{R(\Th)}\sigma_2\mb{}_\bullet \We_J= \ve{R(\Th)}x c_1\mb{}_\bullet \We_J.
\end{eqnarray}
By our assumptions on the action of $\WeD$, the coset $(\ve{R(\Th)}x c_1)_\bullet \We_J$ is a standard parabolic subgroup. It follows $\sigma_1\in\We_J$. By inserting in (\ref{stgg}) and by (iv) we get
\begin{eqnarray*}
   &&   \We_J  =  \ve{R(\Th)}x c_1\mb{}_\bullet \We_J  =  \ve{R(\Th)}x c_1\mb{}_\bullet (\ve{R(J^\infty)}_\bullet\We_J) \\
   &&  =  \ve{R(\Th)}x c_1 \ve{R(J^\infty)}_\bullet\We_J = \ve{R(\Th)}x \ve{R(J^\infty)}c_1\mb{}_\bullet\We_J\\
   &&  =  \ve{R(\Th\cup J^\infty \cup\red{x})}c_1\mb{}_\bullet\We_J = \ve{R(\Th\cup J^\infty \cup\red{x})}_\bullet c_1\We_J.
\end{eqnarray*} 
Applying the idempotent $\ve{R(\Th\cup J^\infty \cup\red{x})}$ on both sides we find
\begin{eqnarray}\label{idbb}
 \ve{R(\Th\cup J^\infty \cup\red{x})}_\bullet\We_J=\ve{R(\Th\cup J^\infty \cup\red{x})}_\bullet c_1 \We_J =\We_J.
\end{eqnarray}
Because in $\WeD$ it holds 
\begin{eqnarray*}
   \We_{\Th\cup J^\infty \cup \red{x}}\ve{R(\Th\cup J^\infty \cup \red{x})}=\ve{R(\Th\cup J^\infty \cup\red{x})},
\end{eqnarray*} 
we get $\We_{\Th\cup J^\infty \cup\red{x}}\subseteq\We_J$ from (\ref{idbb}). It follows $\Th\subseteq J$ and 
$x\in\We_J\cap\mb{}^{\Th\cup\Th^\bot}\We^{J^\infty\cup (J^\infty)^\bot}=\{1\}$. We have shown $\sigma_2=c_1c_2$ with 
$c_1\in(\We_{(J^\infty)^\bot})^{J^0}\subseteq \We_{\Th^\bot}$ and $c_2\in\We_J$. Furthermore, $\sigma_1\in\We_J$. By inserting in (\ref{stgg}) we find
\begin{eqnarray*}
   \We_J=\ve{R(\Th)}\sigma_2\mb{}_\bullet \We_J = \ve{R(\Th)}c_1\mb{}_\bullet \We_J = c_1\ve{R(\Th)}_\bullet \We_J.
\end{eqnarray*} 
Again by our assumptions on the action of $\WeD$, the coset $\ve{R(\Th)}_\bullet \We_J$ is a standard parabolic subgroup. It follows $c_1\in\We_J$.

We have found $\sigma_1\in\We_J$, $\Th\subseteq J$, and $\sigma_2=c_1c_2\in\We_J$, which shows $\hat{\sigma}=\sigma_1\ve{R(\Th)} \sigma_2\in\WeD_J$.
\qed 

Because of the last theorem it is resonable to define:
\begin{defn}\label{DefgoodG} An action of the face monoid $\GD$ on the building $\Omega$, extending the action of the Kac-Moody group $G$ on $\Omega$, is called good if the following two conditions hold:
\begin{itemize}
\item[(1)]$\mb{Stab}_{\widehat{G}}(P_J)=\PD_J$ for all $J\subseteq I$.
\item[(2)] $\GD$ acts order preservingly on $\Omega$.
\end{itemize}
\end{defn}

The face monoid $\GD$ fits very well to the building $\Omega$, which is summarized in:
\begin{cor} The actions of (a) and (b) correspond bijectively by compatibility:
\begin{itemize}
\item[(a)] The actions / order preserving actions / good actions $\bullet$ of $\GD$ on $\Omega$, which extend the action of $G$ on $\Omega$.
\item[(b)] The actions / order preserving actions / good actions $\bullet$ of $\WeD$ on $\mathcal C$, which extend the action of $\We$ on $\mathcal C$, such that $\ve{R(\Th)}_\bullet\sigma\We_J\in {\mathcal C}$ is a standard parabolic for all $\Th\subseteq I$ special, $J\subseteq I$, and $\sigma\in\mb{}^{\Th\cup\Th^\bot}\We^J$.
\end{itemize}
\end{cor}
\Proof It remains to show that an action $\bullet$ of $\GD$ as in (a) is compatible to an action $\bullet$ of $\WeD$ as in (b).
Because of Theorem \ref{AB12} it is sufficient to show $\ND_\bullet{\mathcal A}\subseteq {\mathcal A}$.

Let $\Th\subseteq I$ special, $J\subseteq I$, and $n_\sigma\in N$ projecting to $\sigma\in\mb{}^{\Th\cup\Th^\bot}\We^J$. We first show that $e(R(\Th))_\bullet n_\sigma\PD_J\in\Omega$ is a standard parabolic. Let $e(R(\Th))_\bullet n_\sigma\PD_J= g P_K$. It holds $B=U_\Th U^{\Th\cup\Th^\bot}U_{\Th^\bot}T$ and $\sigma^{-1}U_{\Th^\bot}T\sigma\subseteq B$. By Theorem \ref{FM1} follows for all $u\in U_\Th U^{\Th\cup\Th^\bot}$ and $v\in U_{\Th^\bot}T$:
\begin{eqnarray*}
   uv g P_K &=& uv e(R(\Th))_\bullet n_\sigma\PD_J = u e(R(\Th))v _\bullet n_\sigma\PD_J 
                                          = e(R(\Th))_\bullet n_\sigma (n_\sigma^{-1}v n_\sigma)\PD_J \\
            &=& e(R(\Th))_\bullet n_\sigma\PD_J = gP_K.
\end{eqnarray*}
Therefore, $g^{-1}Bg\subseteq P_K$. It follows $g\in P_K$, compare part 3.2.3 of \cite{Ti}. We have shown $e(R(\Th))_\bullet n_\sigma\PD_J = P_K.$ 

Now $\ND {\mathcal A} \subseteq {\mathcal A}$ follows, because by Theorem \ref{FM1} it holds $e(R(\Th))n = n e(R(\Th))$ for all $n\in N_{\Th\cup\Th^\bot}$, $\Th$ special.\qed 

We now use Theorem \ref{AB12} to show that there are compatible actions of $\GD$ on $\Omega$ for our three actions of $\WeD$ on $\mathcal C$ described before. It is open to determine in particular all good actions of $\GD$ on $\Omega$, which extend the action of $G$ on $\Omega$.

If we apply Theorem \ref{AB12} to the bad action described in Proposition \ref{AC3} we find easily:
\begin{cor}[Bad action] There exists a unique action $\mb{}_\bad$ of $\GD$ on $\Omega$ extending the action of $G$ on $\Omega$, compatible with the bad action $\mb{}_\bad$ of $\WeD$ on $\mathcal C$. 

Explicitely it can be obtained as follows: Let $g_1 e(R(\Th)) g_2\in \GD$ and $h P_J\in \Omega$. 
If $g_2 h P_J\in  P_{\Th\cup \Th^\bot}^-(P_{\Th}\downarrow) =(U^\Th)^-   \We_{\Th^\bot }(P_{\Th}\downarrow)$ decompose $g_2 h$ in the form $g_2h=ab$ with $a\in P_{\Th\cup\Th^\bot}^-$ and $b\in P_J$. Then
\begin{eqnarray*}
  (g_1 e(R(\Th))g_2 )_\bad h P_J =  g_1 p_\Th^-(a) P_J.
\end{eqnarray*} 
If $g_2 h P_J\notin  P_{\Th\cup \Th^\bot}^-(P_{\Th}\downarrow) =(U^\Th)^-   \We_{\Th^\bot }(P_{\Th}\downarrow)$ then
\begin{eqnarray*}
  (g_1 e(R(\Th))g_2 )_\bad h P_J =  P_I.
\end{eqnarray*} 

This action of $\GD$ is in general not order preserving. The parabolic submonoid $\PD_J$ is the stabilizer of $P_J\in\Omega$, $J\subseteq I$.
\end{cor}

If we apply Theorem \ref{AB12} to good action 1 described in Proposition \ref{ACgood1} we find: 
\begin{cor}[Good action 1] There exists a unique action $\mb{}\goodone$ of $\GD$ on $\Omega$ extending the action of $G$ on $\Omega$, compatible with good action 1 $\mb{}\goodone$ of $\WeD$ on $\mathcal C$. It is good.

Explicitely, this action can be obtained as follows: Let $g_1 e(R(\Th)) g_2\in \GD$ and $h P_J\in \Omega$. Decompose $g_2 h$ in the form
\begin{eqnarray*}
  g_2 h= a n_y   c
\end{eqnarray*}
where $a \in P_{\Th \cup \Th^\bot}^-$, $c\in P_J$, and $n_y \in N$ with corresponding projection 
$y\in\mb{}^{\Th\cup\Th^\bot}\We^J$. Then
\begin{eqnarray*}
  (g_1 e(R(\Th)) g_2)\goodone  h P_J = g_1 p_\Th^-(a) P_{\Th\cup J\cup \red{y}}.
\end{eqnarray*} 
\end{cor}

If we apply Theorem \ref{AB12} to good action 2 described in Theorem \ref{AC7} we find: 
\begin{cor}[Good action 2] There exists a unique action $\mb{}_\goodtwo$ of $\GD$ on $\Omega$ extending the action of $G$ on $\Omega$, compatible with good action 2 $\mb{}_\goodtwo$ of $\WeD$ on $\mathcal C$. It is good.

Explicitely, this action can be obtained as follows: Let $g_1 e(R(\Th)) g_2\in \GD$ and $h P_J\in \Omega$. Decompose $g_2 h$ in the form
\begin{eqnarray*}
  g_2 h= a n_x  n_{c_1} c_2
\end{eqnarray*}
where $a \in P_{\Th \cup \Th^\bot}^-$, $c_2\in P_J$, and $n_x,\,n_{c_1}\in N$ with corresponding projections 
$x\in\mb{}^{\Th\cup\Th^\bot}\We^{J^\infty\cup (J^\infty)^\bot}$, $c_1\in (\We_{(J^\infty)^\bot})^{J^0}$ such that in addition 
$x c_1\in \mb{}^{\Th\cup\Th^\bot}\We^J$. 
Set $\Xi:=\Th\cup J^\infty \cup \red{x}$. Then
\begin{eqnarray*}
  (g_1 e(R(\Th)) g_2)_\goodtwo  h P_J = g_1 p_\Th^-(a) P_{\Xi\cup (\Xi^\bot\cap c_1 J^0)}.
\end{eqnarray*}  
\end{cor}

It is instructive to compare the orbit of $B\in\Omega$ for these actions. For simplicity we assume the generalized Cartan matrix $A$ to be indecomposable and of affine or indefinite type. To make our notation plain we denote the action of the Kac-Moody group $G$ on $\Omega$ by a dot ``.''.By applying the face 
monoid $\GD$ to the point $B\in\Omega$ we get for the bad action
\begin{eqnarray*}
 \GD_\bad B&=&G.B\cup P_I,
\end{eqnarray*}
which looks like the flag variety together with an irrelevant point. We get for good action 1
\begin{eqnarray*}
\GD\goodone B&=& \bigcup_{\Th\subseteq I\;special\atop y\in\mb{}^{\Th\cup\Th^\bot}{\mathcal W}}G.P_{\Th\cup\red{y}}.
\end{eqnarray*}
In a subsequent paper we give an algebraic geometric model for a $\GD$-space $\Omega$, where the action of $\GD$ on the building $\Omega$ coincides with good action 1, \cite{M6}. At least this actions occurs quite naturally. We get for good action 2
\begin{eqnarray*}
  \GD_\goodtwo B&=& \bigcup_{\Th\subseteq I\;special} G.P_\Th.
\end{eqnarray*}
%
%
%
%%%%%%%%%%%%%%%%%%%%%%%%%%%%%%%%%%%%%%%%%%%%%%%%%%%%%%%%%%%%%%%%%%%%%%%%%%%%%%%%%%%%%%%%%%%%%%%%%%%%%%%%%%%%%%%%%%%%%%%%%%%%%%%%%%%%%%%%%%%%%
%
%%%%%%%%%%%%%%%%%%%%%%%%%%%%%%%%%%%%%%%%%%%%%%%%%%%%%%%%%%%%%%%%%%%%%%%%%%%%%%%%%%%%%%%%%%%%%%%%%%%%%%%%%%%%%%%%%%%%%%%%%%%%%%%%%%%%%%%%%%%%%%%%%
%
%

%
\end{document}